\documentclass[a4paper,12pt,notitlepage,english,onecolumn,oneside,fleqn]{article}
\author{Israel Moreno-Mej\'{\i}a, Dan Silva-L\'opez}
\usepackage{pgf}                
\usepgflibrary{arrows}
\usepackage{tikz}                
\usetikzlibrary{shapes,arrows}
\usetikzlibrary{arrows,decorations,backgrounds}
\usetikzlibrary{babel}
\usepackage{mathrsfs}
\usepackage[cal=boondox]{mathalfa}
\usepackage{eufrak}
\usepackage{latexsym}
\usepackage{pst-all}
\usepackage{amscd}
\usepackage{amsfonts}
\usepackage{alphalph}
\usepackage{amsmath}
\usepackage{ifsym}
\usepackage{dsfont}
\usepackage[toc]{appendix}
\usepackage{euscript}
\usepackage{fancyhdr}
\usepackage{makeidx}
\usepackage{amsthm}
\usepackage{bigints}
\usepackage[latin1]{inputenc}
\usepackage{amssymb, amsmath, amsbsy} 
\usepackage{cancel}
\definecolor{light-gray}{gray}{0.95}
\definecolor{orange1}{cmyk}{0,0.5,1,0}
\definecolor{orange2}{rgb}{1,0.5,0}
\definecolor{yellow}{cmyk}{0,0,0.15,0.05}
\makeindex
\makeglossary
\newtheorem{d1}{Definici\'on}[section]

\newtheorem{t1}[d1]{Theorem}

\newtheorem{l1}[d1]{Lemma}

\def\<{\langle}
\def\>{\rangle}
\newcommand{\nothing}[1]{}
\newcommand{\mathftnt}[1]{\mbox{\footnotesize$#1$}}

\newcommand{\mathtiny}[1]{\mbox{\tiny$#1$}}
\newcommand{\mathsmall}[1]{\mbox{\small$#1$}}
\newcommand{\LinL}[1]{\overline{{#1} \hspace{1pt}}}
\usepackage{hyperref} 
\numberwithin{equation}{section}
\newdimen\arrowsize
\pgfarrowsdeclare{inyectivo1}{inyectivo1}{}
{
\pgfsetdash{}{0pt} 
\pgfsetroundjoin
\pgfsetroundcap
\pgfpathmoveto{\pgfpoint{-10pt}{5pt}}
\pgfpatharc{90}{180}{11pt and 4pt}
\pgfpatharc{180}{197}{17pt and 3pt}
\pgfusepathqstroke
\pgfpathmoveto{\pgfpoint{-4pt}{4pt}}
\pgfpatharc{180}{270}{4pt}
\pgfpatharc{90}{180}{4pt}
\pgfusepathqstroke
}
\newdimen\arrowsize
\pgfarrowsdeclare{inyectivo2}{inyectivo2}{}
{
\pgfsetdash{}{0pt} 
\pgfsetroundjoin
\pgfsetroundcap
\pgfpathmoveto{\pgfpoint{0pt}{.1pt}}
\pgfpatharc{-30}{70}{9pt and 3.3pt}
\pgfpatharc{75}{110}{21pt and 7pt}
\pgfusepathqstroke
\pgfusepathqstroke
}
\newdimen\arrowsize
\pgfarrowsdeclare{inyectivo3}{inyectivo3}{}
{
\pgfsetdash{}{0pt} 
\pgfsetroundjoin
\pgfsetroundcap
\pgfpathmoveto{\pgfpoint{0pt}{.1pt}}
\pgfpatharc{30}{-70}{9pt and 3.3pt}
\pgfpatharc{75}{110}{21pt and 7pt}
\pgfusepathqstroke
\pgfusepathqstroke
}
\newdimen\arrowsize
\pgfarrowsdeclare{arcs}{arcs}{}
{
\pgfsetdash{}{0pt} 
\pgfsetroundjoin 
\pgfsetroundcap 
\pgfpathmoveto{\pgfpoint{-4pt}{4pt}}
\pgfpatharc{180}{270}{4pt}
\pgfpatharc{90}{180}{4pt}
\pgfusepathqstroke
}
\pgfarrowsdeclare{arcs'}{arcs'}
{
\arrowsize=0.2pt
\advance\arrowsize by .6\pgflinewidth
\pgfarrowsleftextend{-6\arrowsize-.6\pgflinewidth}
\pgfarrowsrightextend{.6\pgflinewidth}
}
{
\arrowsize=0.2pt
\advance\arrowsize by .6\pgflinewidth
\pgfsetdash{}{0pt} 
\pgfsetroundjoin
\pgfsetroundcap
\pgfpathmoveto{\pgfpoint{-6\arrowsize}{6\arrowsize}}
\pgfpatharc{180}{270}{6\arrowsize}
\pgfusepathqstroke
\pgfpathmoveto{\pgfpointorigin}
\pgfpatharc{90}{180}{6\arrowsize}
\pgfusepathqstroke
}
\newdimen\arrowsize
\pgfarrowsdeclare{supra}{supra}{}
{
\pgfsetdash{}{0pt} 
\pgfsetroundjoin
\pgfsetroundcap
\pgfpathmoveto{\pgfpoint{-4pt}{4pt}}
\pgfpatharc{180}{270}{4pt}
\pgfpatharc{90}{180}{4pt}
\pgfpathmoveto{\pgfpoint{-7pt}{4pt}}
\pgfpatharc{180}{270}{4pt}
\pgfpatharc{90}{180}{4pt}
\pgfusepathqstroke
}
\pgfarrowsdeclarereversed{in2 reversed}{in2 reversed}{inyectivo2}{inyectivo2}
\title{The Verlinde traces for $SU_{X}(2,\xi)$ and blow-ups}
\begin{document}
\maketitle
\begin{abstract}\footnote{2000 MSC 14F25,14H60,32L05,14H37,(14F05,14D20,14C17,14H45,20C15)}
Given a compact Riemann surface $X$ of genus at least $2$ with automorphism group $G$
we provide formulae that enable us to compute traces of automorphisms of  $X$  on the space of global sections of  $G$-linearized line bundles defined on certain blow-ups of projective spaces along the curve $X$. The method is an adaptation of one used by Thaddeus to compute the dimensions of those spaces. In particular we can compute the traces of automorphisms of  $X$ on the Verlinde spaces corresponding to the moduli space  $SU_{X}(2,\xi)$ when $\xi$ is a line bundle $G$-linearized  of suitable degree.  
\end{abstract}
\section{Introduction}

Let $X$ be a complex, irreducible, smooth, projective curve of genus at least $2$ and 
automorphism group $G=Aut(X)$.  Let $\xi$ be a $G-$linearized line bundle over $X$.

By the Verlinde traces we refer to the traces of 
 automorphisms of $X$ on the space $H^0(SU_{X}(r,\xi),\mathcal{O}(n\Theta))$,  
where $SU_{X}(r,\xi)$ is the moduli space of semi-stable rank r vector bundles with determinant $\xi$  and 
where $\mathcal{O}(\Theta)$ is the determinantal line bundle of $SU_{X}(r,\xi)$. In this work we address the problem of 
computing the Verlinde traces for the case $r=2$ and we give a way to compute them which is alternative to the approach of 
J. E. Andersen for case of the trace $\tau_{\mathcal{G}}^{(k)}(f)$   in \cite{MR3181488} pg. 3 ($\mathcal{G}= SU(2)$ in our case).     

Our method and result will be explained in the  paragraphs below. Before that, we would like to mention a few things related to this problem. As the reader may be aware the case of the Verlinde traces for the identity of $G$ is already solved. A   formula for
 $dim$ $H^0(SU_{X}(r,\xi),\mathcal{O}(n\Theta))$ was conjectured by E. Verlinde \cite{MR954762} and there are many proofs of 
it(see for instance the following works concerning  this case \cite{MR1231965},\cite{MR1231966},\cite{MR1273268}, \cite{MR1215963}, 
\cite{MR1244913}, \cite{MR1245399}, \cite{MR1360519}, \cite{MR1257326}, \cite{MR1289330}, \cite{MR1048605} and \cite{MR1133264}). 
The action of non-trivial automorphism groups $G$ on the Verlinde spaces had already been considered in the work of Dolgachev 
(\cite{MR1653063} in his proof of Cor. 6.3), 
 in the work of the first author  \cite{MR2125533} (see Tables 1--4 there) where some Verlinde traces are 
computed for the cases $\xi=\mathcal{O}_{X}$, $n=1$ and arbitrary rank $r$ by computing them on 
$H^0(J_{X}^{g_{X}-1},\mathcal{O}(r\Theta_{J_{X}^{g_{X}-1}}))^{*}\cong H^0(SU_{X}(r,\xi),\mathcal{O}(\Theta))$ and later in the work of  Andersen \cite{MR3181488} the Verlinde traces correspond to the traces  $\tau_{\mathcal{G}}^{(k)}(f)$   which we will briefly  refer to  at the end of this introduction.

As there are automorphisms of $SU_{X}(r,\xi)$ that are not induced by $G$ (for a description $Aut(SU_{X}(r,\xi))$ and 
related results see \cite{MR1336336},\cite{MR2103475}, \cite{MR3035121}) we should also mention that the action of torsion 
elements of the Jacobian of $X$ acting on the Verlinde spaces of some moduli spaces of vector bundles  had also been 
considered in the works  \cite{MR424819}, \cite{MR2795755} and \cite{MR1490854}. Explicit formulae for the corresponding 
Verlinde traces are provided in the latter two. \\
Coming back to our problem, in the rank 2 case,  we followed the method used by Thaddeus  to derive the Verlinde formula in \cite{MR1273268}. We shall see that his method can be extended to compute the Verlinde traces 
(see formula (\ref{eq:3c2p1}) and Section \ref{TPoThaddeus})
by just replacing the use of the Riemann-Roch Theorem  
for the use of the Atiyah-Singer Holomorphic Lefschetz Theorem 
(see Section \ref{sec:THL}) and in this work  
we derive some formulae (Theorem \ref {t:chhBimn} ) required to apply the Holomorphic Lefschetz Theorem in Thaddeus' method.

Let $K_{X}$ be the canonical line bundle of $X$. Suppose that $K_{X}\xi $ is very ample. Let
 $X\hookrightarrow\mathds{P}^{N}$
be the embedding defined by the complete linear system $\mid K_{X}\xi\mid$.
Let $\pi:\widetilde{\mathds{P}^{N}_{X}}\mapsto \mathds{P}^{N}$ be the blow-up of $\mathds{P}^{N}$ with center $X$ and let $E$ be the corresponding exceptional divisor. The Picard group of $\widetilde{\mathds{P}^{N}_{X}}$ is generated by  $\mathcal{O}(E)$ and the hyperplane line bundle $\mathcal{O}(H)$. For integers $m,n$ let 
$\mathcal{O}_{1}(m,n)= \mathcal{O}((m+n)H-nE)$ and let 
$V_{m,n}=H^{0}(\widetilde{\mathds{P}^{N}_{X}}, \mathcal{O}_{1}(m,n))$.
 Let $\xi$ have degree $d$. Thaddeus shows that
for $d>2g_{X}-2$ there is a natural isomorphism
\begin{equation}\label{eq:espaciosVerlinde}
H^0(SU_{X}(2,\xi),\mathcal{O}(k\Theta))\cong V_{k,k(d/2-1)}.
\end{equation}

Under some mild conditions on $m, n $ he  finds a formula for the dimension of 
 $V_{m,n}$ (see Theorem \ref{dimvmn} below). The cases not covered by those conditions can be dealt with easily.

Now, when we assume that $\xi $ is $G$-linearized   it induces   an action of $G$ on $\mathds{P}^{N}$ and on the blow-up 
$ \widetilde{\mathds{P}^{N}_{X}}$ such that the embedding $X\hookrightarrow\mathds{P}^{N}$  and the blow-up map $\pi$ are $G$-equivariant.  
The line bundles  $\mathcal{O}_{1}(m,n)$ can be equipped with a linearization induced by that of $\xi$. 
That is because  $\mathcal{O}(H)= \pi^{*}{\mathcal{O}_{\mathds{P}^{N}}(1)}$ comes equipped with a linearization induced by that 
of $\xi$; also since $E$ is a $G$-invariant divisor $\mathcal{O}(E)$  admits a linearization of $G$ with trivial action on $H^0( \widetilde{\mathds{P}^{N}_{X}},\mathcal{O}(E))$,
and since $h^0( \widetilde{\mathds{P}^{N}_{X}},\mathcal{O}(E))=1$ we can see that this 
linearization is unique(with trivial action on $H^0( \widetilde{\mathds{P}^{N}_{X}},\mathcal{O}(E))$)  because 
there is a bijection between the $G$-linearizations of $\mathcal{O}(E)$ and the 
$G$-invariant divisors linearly equivalent to E by $G$-invariant rational functions (see \cite{MR1653063} Prop. 2.1 and \cite{MR2004511} Ex. 7.4).
 So one can consider the problem of computing the traces of elements of $G$ on the spaces $V_{m,n}$ (Thaddeus traces). 
 As we will see in Section \ref{stf}, the formula for $dim V_{m,n}$ is a linear combination of Euler Characteristics of sheaves  $B_{i,m,n}$ defined over symmetric products $S^{i}X$ of the curve (see  Theorem \ref{dimvmn} below) and these sheaves are naturally  $G$-linearized. By tracking back the proof of  Theorem \ref{dimvmn} one notice that the homomorphisms between the cohomology groups involved  are $G$-equivariant and that the trace of an element  $h\in G$ on $V_{m,n}$ is in fact obtained by replacing the Euler characteristics of the sheaves $B_{i,m,n}$ by their corresponding Lefschetz numbers $N_{i}(h)=L(h,B_{i,m,n})$ (see formula (\ref{eq:3c2p1}) below and Section \ref{TPoThaddeus} for its proof).
These Lefschetz numbers $L(h,B_{i,m,n})$ can be computed using the Holomorphic Lefschetz Theorem (see  Section \ref{sec:THL}, Theorem \ref{teo1} and formula (\ref{eq:thlsum})) because the sheaves $B_{i,m,n}$ are defined on smooth varieties. To do that we need to know for each component $Z$ of the fixed point set the following data:
the generalized Chern Character $\textnormal{ch}_{h}(i^{*}_{Z}B_{i,m,n})$ (see (\ref{eq:defchg})), 
the stable characteristic classes $U(N_{Z/S^{i}X}(\nu^{j}))$ of the normal bundle $N_{Z/S^{i}X}$
(see (\ref{eq:ec35})), the Todd class  $\textnormal{Td}(T_{Z})$ of the tangent bundle $T_{Z}$ and $\textnormal{det}(Id-h_{\mid}{N^{\vee}_{Z/S^{i}X}})$ (see Theorem \ref{teo1}).
In this way by equation (\ref{eq:thlsum}) the solution of the problem of computing $N_{i}(h)=L(h,B_{i,m,n})$ is 
reduced to the calculation of
the following intersection numbers, namely, the contributions
\begin{equation}\label{eq:ec73}C_{i,Z}(h):=\bigints_{Z}
\frac{\textnormal{ch}_{h}(i^{*}_{Z}B_{i,m,n})[\prod_{j=1}^{o(h)-1}\textnormal{U}
(N_{Z/{S^{i}X}}(\nu^{j}))]\textnormal{Td}(T_{Z})}
{\textnormal{det}(Id-h_{\mid}{N^{\vee}_{Z/{S^{i}X}}})}
\end{equation}
of each component $Z\subseteq {(S^{i}X)}^{h}$ of the fixed point set 
of $h.$ 
\normalsize
 As we pointed out earlier, the main goal of this paper (Theorem \ref {t:chhBimn}) is the calculation of the generalized  Chern Character $\textnormal{ch}_{h}(i^{*}_{Z}B_{i,m,n})$. The other data required to apply (\ref{eq:ec73}) 
have been dealt with in the works \cite{MR2125533}(Prop. 3.2) and \cite{MR2254543} (Theorem 2.3). However, in Section \ref{SCC} we present a generalization of the formula for the stable characteristic classes $U(N_{Z/S^{i}X}(\nu^{j}))$ that was given in \cite{MR2254543}. The Theorem \ref{chfw+w-} is required for the proof of Theorem \ref{t:chhBimn} and part of its proof is modelled on the proof  of Proposition 2.1 in \cite{MR2254543}. At the end of the paper we illustrate  the use of formula (\ref{eq:3c2p1}), when $X$ is a hyperelliptic curve of genus 2, by computing the Verlinde traces corresponding to the hyperelliptic involution. 

The problem of computing  Verlinde traces  may be considered interesting in its own right and following the work of Thaddeus on moduli of pairs there have been various constructions 
of series flips using moduli of pairs, triplets, that one should expect to be able to use to compute Verlinde traces on  moduli spaces other than $SU_{X}(2,\xi)$.
Some of our motivations are to determine the isotypic decomposition of the Verlinde spaces
 $H^0(SU_{X}(r,\xi),\mathcal{O}(n\Theta))$ to study the varieties defined as the zero locus of the invariant submodules in, for instance,
 $\mathds{P}H^0(SU_{X}(r,\xi),\mathcal{O}(2\Theta))^{\vee}$; to study the Cox Ring of Blow-ups( in this case with the Thaddeus traces).
Also, as it can be seen from the work of Andersen,   Verlinde traces are used to compute invariants of 
3-Manifolds, namely, the Witten-Reshetikhin-Turaev invariants of finite order mapping tori
 $Z_{G}^{(k)}(\sum f)=Det(f)^{-\frac{1}{2}\zeta}\tau_{\mathcal{G}}^{(k)}(f)$.
Andersen considers moduli spaces of semi-stable $\mathcal{G}$-bundles on curves and by applying directly the Lefschetz-Riemann-Roch Theorem  for singular varieties (see \cite{MR549774}) 
to the determinantal line bundle of the moduli spaces he obtained an expression for the traces $\tau_{\mathcal{G}}^{(k)}(f)$ because the higher cohomologies of the determinantal line bundle vanish. The formula is given by the equation 8.1 or Theorem 1.3 in \cite{MR3181488} up to the correction
factor $Det(f)^{-\frac{1}{2}\zeta}$,
 and depends on data coming from the components of 
the fixed point set. The moduli spaces $SU_{X}(r,\mathcal{O}_{X})$ correspond to the case $\mathcal{G}=SU(r)$
and a description of the components of fixed point set for this case is given in Theorem 6.10 of \cite{MR3181488}
(see also Theorem 3.4 of \cite{MR2204258}). 
The contribution of the smooth components to the traces $\tau_{\mathcal{G}}^{(k)}(f)$ is also determined in his work. In the rank 2 case, the advantage of our method is that one does not have to deal with singular components of fixed points and one can compute the Verlinde traces if  one knows the action of the automorphism on the tangent spaces of the fixed points in the curve.

\vspace{0.5cm}
\noindent Some results of this work  are based on results of the Ph.D Thesis of the second author \cite{silva}.

\newpage

\section{Notation}
Let $p$ be a positive integer and let $\nu=\exp(2i\pi/p)$.
Given a finite cyclic group $H=\< h \> $ of order  $p$ acting trivially on a variety $Z$ and given  an $H$-linearized vector bundle $F$ on $Z$ there is a decomposition into eigenbundles
\[F=\bigoplus_{j=0}^{p-1} F(\nu^j),\]
that is, $F(\nu^j)$ is the sub-bundle of $F$ where the action of $h$  on the fibres is multiplication by
$\nu^j$. Most of the times we will say that $F$ is $h$-linearized.

For a divisor  $D$ on a variety $W$ we write $\mathcal{O}_{W}(D)$ (or just $\mathcal{O}(D)$) to denote corresponding line-bundle of $D$. If $F$ is a sheaf on  $W$ we usually write  $F(D)$ rather than $F\otimes \mathcal{O}(D)$. Usually $F^{n}=F^{\otimes n}$, also if  $K$ is another sheaf some times we could write $KF= K\otimes F$.  
If $\mathscr{\xi}$ is a sheaf on $W$ and $\xi'$ is a sheaf on 
$X$ we use the notation
$\xi \boxtimes \xi'=\pi^{*}_{W}\xi \otimes \pi^{*}_{X}\xi'$ on the product 
$W\times X$, where $\pi_{W}$ and 
$\pi_{X}$ are the natural projections and if $G$ is a sheaf on $W\times X$ we sometimes write $\xi G=\pi_{W}^{*}\xi\otimes G$ (or 
$\xi (D)$ instead of $\xi G$
if $G=\mathcal{O}_{W\times X}(D)$ with $D\subset W\times X$
is a divisor).  

If $i:\Delta\hookrightarrow W$ is the inclusion of a subvariety of $W$ then $\mathcal{O}_{\Delta} F=\mathcal{O}_{\Delta}\otimes F$ is the pull-back
$i^{*}F$ for any sheaf $F$ on $W$ and if $\Delta$ is a divisor of $W$ then we denote
$\mathcal{O}_{\Delta}(\Delta)=i^{*}\mathcal{O}_{W}({\Delta})$. If $F$ is locally free then $F^{\vee}$ is the dual sheaf and we also denote by $F$ the vector bundle  whose sheaves of sections is $F$. If $x\in W$ then  $F_{x}$ is the fibre of $F$ at the point $x$. For a vector bundle $F$ over $X$ we denote by $\pi:\mathds{P}(F)\rightarrow X $ the projective bundle of one dimensional subspaces of $F$ thus $\mathds{P}(F)=Proj_{\mathcal{O}_{X}}(Sym(F^{\vee}))$ where $Sym(E)=\oplus_{m} S^{m}E$ is the symmetric algebra of a locally free sheaf $E$ over $X$. The exponential function $exp(x)=\sum_{n=0}^{\infty}x^{n}/n!$ is also be represented by $e^{x}$ and in some formulae we mixed both notations.

\section{The Holomorphic Lefschetz Theorem}\label{sec:THL}
Consider a finite cyclic group $H=\< h \> $ of order  $p\geq1$.\\
 Suppose that $H$ is acting trivially on a variety $Z$ and also consider an $h$-linearized vector bundle $F$ on $Z$. 
The generalized  Chern character of $F$ is given by
\begin{equation}\label{eq:defchg}
ch_{h}(F)=\sum_{j=0}^{p-1}\nu^{j}ch[F(\nu^{j})],
\end{equation}
 where $ch[F(\nu^{j})]$ is the Chern character of the eigenbundle $F(\nu^{j})$.
 For each vector bundle  $F(\nu^{j})$  with $\nu^{j}\not=1$ define the stable characteristic class
  $U(\!F(\nu^{j})\!)$ as
\begin{equation}\label{eq:ec35}
U(F(\nu^{j}))=\prod_{i=1}^{r_{j}}\left(\frac{1-\frac{e^{-x_{i}}}{\nu^{j}}}{1-\frac{1}{\nu^{j}}}\right) ^{-1},
\end{equation}
where $r_{j}$ is the rank $F(\nu^{j})$ and the $x_{i}$'s are the Chern roots of $F(\nu^{j})$.
\\
 We now make some observations about  fixed point sets of  $h$ acting on  non singular Varieties (see pg. 537 in \cite{MR236951} and Lemma 4.1  in \cite{MR263834}).
Let $X$ be an irreducible, non-singular complex algebraic variety on which $h$ acts as an automorphism of $X$. Let $X^{h}$ denote the set of fixed points of $h$ in $X$.  One has that    $X^{h}$ is non-singular and if $x\in X^{h}$ then  as an analytic variety $X^{h}$ is (in a neighbourhood of $x$) the image, under a suitable equivariant exponential map,  of the eigenspace $T_{X,x}(\nu^{0})$ of the tangent space $T_{X,x}$ to $X$ at $x$. As $T_{X^{h},x}=T_{X,x}(\nu^{0})$ one has that the linear transformation $h_{\mid} N_{X^{h}/X,x}$ induced by $h$ on the normal space $ N_{X^{h}/X,x}$ to $X^{h}$ at $x$ has no eigenvalue $\nu^0$ and therefore
\begin{equation}
\det(Id-h\mid N_{X^{h}/X,x})\not=0.
\end{equation}
 One also notice from what was  said above that  $X^{h}$ is locally irreducible so a point $x\in X^{h}$ belongs to a unique irreducible component of $X^{h}$ and so
 \begin{equation}\label{eq:cpdirrc}
 X^{h}=\coprod_{Z\in CX^{h}} Z,
 \end{equation}
 where $CX^{h}$ is the set of irreducible components of $X^{h}$.

\vspace{0.5cm}
Let $i_{X^{h}}:X^{h}\hookrightarrow X$ the inclusion $X^{h}\subseteq X$.
Let $E$ be an $h$-linearized vector bundle on $X$.
The Lefschetz number $L(h,E)$ (sometimes written $L(h,X,E)$) of $h$ on $E$ 
is 
\[ L(h,E)= \sum_{i}(-1)^{i}trz\ h|H^{i}(X,E)  \]
and can be computed by the \emph{Holomorphic Lefschetz Theorem}:
\begin{t1}\label{teo1}%
\textbf{(\small Atiyah-Singer,
\textit{Theorem 4.6, pg.566, \cite{MR236952}})} 
\\
\normalsize
\footnotesize
\begin{equation}\label{eq:ec73ppr}
\begin{split}
L(h,E)&=
\bigint_{X^{h}}\frac{\textnormal{ch}_{h}(i^{*}_{X^{h}}E)[\prod_{j=1}^{o(h)-1}\textnormal{U}
(N_{X^{h}/X}(\nu^{j}))]\textnormal{Td}(T_{X^{h}})}
{\textnormal{det}(Id-h_{\mid}{N^{\vee}_{X^{h}/X}})}\\
\end{split}\\
\end{equation}

\normalsize
\noindent where ${\textnormal{det}(Id-h_{\mid}{N^{\vee}_{X^{h}/X}})}$  assigns to the component of $x\in X^{h}$ the value of
${\textnormal{det}(Id-h_{\mid}{N^{\vee}_{X^{h}/X,x}})}.$ 
 Also,   $U(N_{X^{h}/X}(\nu^{j}))$ 
is the stable characteristic class of $N_{X^{h}/X}(\nu^{j})$, \textnormal{ch}$_{h}(i^{*}_{X^{h}}E)$ is the generalized Chern character of the pull-back $i^{*}_{X^{h}}E$ of $E$ to $X^{h}$ and $\textnormal{Td}(T_{X^{h}})$ is the Todd class of the tangent sheaf  $T_{X^{h}}$ of $X^{h}$.          
\end{t1}
On the right-hand side of (\ref{eq:ec73ppr}) one is evaluating a cohomology class    
$u\in H^{*}(X^{h},\mathbb{C}) $
(here   $u$  is the integrand).
By (\ref{eq:cpdirrc}) on has that
\[H^{*}(X^{h},\mathbb{C})=\bigoplus_{Z\in CX^{h}} H^{*}(Z,\mathbb{C})\]
and then the evaluation of $u$ is the sum of the evaluations of the components $u_{Z}\in H^{*}(Z,\mathbb{C})$  of $u$, that is,

\begin{equation}\label{eq:thlsum}
L(h,E)=\sum_{Z\in \textnormal{CX}^{h}}\bigints_{Z}
\frac{\textnormal{ch}_{h}(i^{*}_{Z}E)[\prod_{j=1}^{o(h)-1}\textnormal{U}
(N_{Z/X}(\nu^{j}))]\textnormal{Td}(T_{Z})}
{\textnormal{det}(Id-h_{\mid}{N^{\vee}_{Z/X}})},
\end{equation}
A particular case of Theorem \ref{teo1} is the \emph{Atiyah-Bott formula} also known as the \emph{Woods 
Hole Fixed Point Theorem} which is obtained when $X^{h}$ is a finite set:

\begin{equation}\label{eq:z274ppr}
L(h,E)=\sum_{x\in X \mbox{$^{h}$}}\frac{trz \ h|E_{x}}{det(Id-h|T^{\vee}_{X,x})},
\end{equation}
where $trz \ h|E_{x}$ is the trace of $h|E_{x}$(
\textit{see  \cite{MR3848672} pgs. 631 and 121}).

\section{The trace formula}\label{stf}

From now on what follows $X$ will denote the curve considered in the
introduction, that is, a complex, irreducible, smooth curve of genus $g\geq 2$ embedded in $\mathds{P}^{N}=\mathds{P}H^{1}(\xi^{-1})$ via the linear  system $|K_{X}\xi|$, where $\xi$ is a line bundle on $X$ of degree $d$. 
\\
Let $\pi_{S^{i}X} : X\times S^{i}X\mapsto S^{i}X$ and 
$\pi_{X}:X\times S^{i}X\mapsto X$ be the natural projections. Let $\Delta_{i}\subset X\times S^{i}X$ be the universal divisor and let $j'$ denote its inclusion into $X\times S^{i}X$. Consider the \emph{Thaddeus bundles}
\begin{equation}\label{eq:wimen}
\begin{array}{lll}
W^{-}_{i}&=& (R^0_{\pi_{S^{i}X}})_{*}\mathcal{O}_{\Delta_{i}}\xi(-\Delta_{i})\\
&=&(R^0_{\pi_{S^{i}X}})_{*}\{\{j'_{*}\mathcal{O}_{\Delta_{i}}\}\otimes \pi_{X}^{*}(\xi)\otimes\mathcal{O}_{X\times S^{i}X}(-\Delta_{i})\}\\
&=&(R^0_{\pi_{S^{i}X}})_{*}\{\{j'_{*}\mathcal{O}_{\Delta_{i}}\}\otimes (\xi\boxtimes\mathcal{O}_{S^{i}X})(-\Delta_{i})\}
\end{array} 
\end{equation}
and
\begin{equation}\label{eq:wimas}
\begin{array}{lll}
W^{+}_{i}&=&(R^1_{\pi_{S^{i}X}})_{*}\xi^{-1}(2\Delta_{i}) \\
&=&(R^1_{\pi_{S^{i}X}})_{*}\{\pi_{X}^{*}(\xi^{-1})\otimes\mathcal{O}_{X\times S^{i}X}(2\Delta_{i})\}\\
&=&(R^1_{\pi_{S^{i}X}})_{*}\{(\xi^{-1} \boxtimes\mathcal{O}_{S^{i}X})(2\Delta_{i})\}.
\end{array} 
\end{equation}
\noindent
One has that $W^{-}_{i}$ is a vector bundle of rank $i$ 
and that if $-d+2i < 0$ then $W^{+}_{i}$ is a vector bundle of rank $d+g-1-2i$. 
\\
Define
\begin{equation}\label{eq:z56}
	L_{i}=det^{-1}[(\pi_{S^{i}X})_{!}
	\{\xi\mathcal{O}_{X\times S^{i}X}(-\Delta_{i})\}]\otimes det^{-1}[(\pi_{S^{i}X})_{!}
	\{\mathcal{O}_{X\times S^{i}X}(\Delta_{i})\}],
	\end{equation}
	where $(\pi_{S^{i}X})_{!}:K(X\times S^{i}X)\rightarrow K(S^{i}X)$ is the direct image homomorphism between  Grothendieck groups   given by 
	\begin{equation}
	(\pi_{S^{i}X})_{!}(\mathscr{F})
	=\sum_{i}(-1)^{i}{R^{i}_{\pi_{S^{i}X}*}}
	(\mathscr{F})
	\end{equation}
	for a coherent sheaf $\mathscr{F}$ over 
	$X\times S^{i}X$ (See \cite{MR0463157} pg. 436). Also, $\det:K(S^{i}X)\rightarrow Pic(S^{i}X)$ is the group homomorphism 
	between the additive group $K(S^{i}X)$ and the multiplicative group $Pic(S^{i}X)$ given by
	\begin{equation}
	det(\sum_{j=1}^{m}a_{j}\mathcal{F}_{j})=\bigotimes_{j=1}^{m}det(\mathcal{F}_{j})^{a_{j}},
	\end{equation}
	where $a_{1},\cdots, a_{m}\in\mathbb{Z}$, the $\mathcal{F}_{j}'s$ are coherent sheaves and $\det(\mathcal{F}_{j})$ is the usual determinant of coherent sheaves defined by means of a locally free resolution of $\mathcal{F}_{j}$.\\
	Let $U_{i}\rightarrow S^{i}X$ be the bundle 
	\begin{equation}\label{eq:z55}
	U_{i}=W_{i}^{-}\oplus (W_{i}^{+})^{\vee}.
	\end{equation}
Let $q_{i}=n-(i-1)m$. For $i>0$ consider the Euler characteristic
	\begin{equation}\label{eq:z57}
	N_{i}=\chi(S^{i}X,B_{i,m,n}),
	\end{equation}
	where
	\begin{equation}\label{eq:z58}
	B_{i,m,n}=L_{i}^{m}\otimes \wedge^{i}W_{i}^{-}\otimes S^{q_{i}-i}U_{i}.
	\end{equation}
	 For $i=0$ we define $B_{0,m,n}=S^{m+n}H^{0}(X,K_{X}\xi)$ and
	 $N_{0}=dim B_{0,m,n}$. It is a convention that $B_{i,m,n}=0$ when $q_{i}-i<0$. 
	
	 \begin{t1}\label{dimvmn}(See (6.9) in\cite{MR1273268}) Let $m,n\geq 0$ and suppose that  \\
	 $m(d-2)-2n   >  -d+2g-2$. Then 
\begin{equation*}
	dim V_{m,n}=\sum_{i=0}^{\infty}(-1)^{i}N_{i}=\sum_{i=0}^{w}(-1)^{i}N_{i},
	\end{equation*}
	where $w=[(d-1)/2]$.
 \end{t1}
 \noindent Under our hypothesis (namely, those of Theorem \ref{dimvmn}  together with the one  that $X$ has 
 automorphism group $G$ and that $\xi$ is 
 $G$-linearized so that line bundles $\mathcal{O}_{1}(m,n)$ 
 are $G$-linearized) one has that for any $h\in G$ (see Section \ref{TPoThaddeus})
 \begin{equation}\label{eq:3c2p1}
	Trace( h_{\mid_{ V_{m,n}}} ) =\sum_{i=0}^{\infty}(-1)^{i}N_{i}(h)=\sum_{i=0}^{w}(-1)^{i}N_{i}(h),
	\end{equation}
 where $N_{i}(h)$ stands for the Lefschetz number
 \begin{equation}\label{eq:3c2p2}
 L(h,B_{i,m,n} )=\sum_{j=0}^{i} (-1)^{j} Trace(h_{\mid_{}}{{ H^j(S^{i}X,B_{i,m,n}) }}),
 \end{equation}
 and for $i=0$, 
\begin{equation}\label{eq:3c2p2i0}
N_{0}(h)=Trace(h_{\mid_{}}B_{0,m,n})=\underset{t^{m+n}}{coef}\left[\frac{1}{\det(I-t\cdot h_{\mid_{}} H^{0}(X,K_{X}\xi))}\right].
\end{equation} 
The far right-hand side of (\ref{eq:3c2p2i0}) follows from the proof of  \emph{Molien's Theorem 1.10} in  \cite{MR2004218}.
So the Lefschetz numbers $N_{i}(h)$, $i>0$, can be computed by means of the Holomorphic Lefschetz Theorem if one (see formula (\ref{eq:thlsum}))  can determine the contribution 
$$C_{i,Z}(h)=\bigints_{Z}
\frac{\textnormal{ch}_{h}(i^{*}_{Z}B_{i,m,n})[\prod_{j=1}^{o(h)-1}\textnormal{U}
(N_{Z/{S^{i}X}}(\nu^{j}))]\textnormal{Td}(T_{Z})}
{\textnormal{det}(Id-h_{\mid}{N^{\vee}_{Z/{S^{i}X}}})}$$
of each component $Z\subseteq {(S^{i}X)}^{h}$ of the fixed point set.
 It will be explained in Section \ref{sec:TchC} that the  components $Z$ of ${(S^{i}X)}^{h}$ are parametrized by  certain set of $h$-invariant divisors $D$ on the curve $X$ so in the subsequent sections we will write $Z_{D}$ rather than just $Z$ because although components associated to distinct divisors may be isomorphic the data associated to them may not be the same.
  So far, the main obstruction to compute
$C_{i,Z}(h)$ is the generalized Chern character $\textnormal{ch}_{h}(i^{*}_{Z}B_{i,m,n})$ and it will be determined in Section \ref{TgCcfBimn}. 
In Section \ref{SCC},  Theorem \ref{tt3} we present a new formula for the stable characteristic classes
 $\textnormal{U}(N_{Z/{S^{i}X}}(\nu^{j}))$ and in (\ref{aival}) one for ${\textnormal{det}(Id-h_{\mid}{N^{\vee}_{Z/{S^{i}X}}})}.$ The Todd class of a component is given in formula (\ref{eq:97}).

\section{The proof of Thaddeus}\label{TPoThaddeus}
Here we shall derive formula (\ref{eq:3c2p1}) so we are assuming that a finite group $G$ is acting on the curve $X$ and that $\xi$ is a $G$-linearized line bundle of degree $d$. Recall also that  $K_{X}\xi $ is very ample and $X$ is embedded into $\mathds{P}H^{1}(X,\xi^{-1})$  by the complete linear system $\mid K_{X}\xi\mid$. We should remark that $\xi$ plays the role of the line bundle $\Lambda$ considered in \cite{MR1273268}.
Thaddeus' proof of Theorem \ref{dimvmn} uses 4 results ( (6.2),(6.6),(6.7) and (6.8) in \cite{MR1273268} ) which we will state here in terms of Lefschetz numbers ( see Lemmas \ref{6.2thadppr}--\ref{6.8thadppr} below), their proofs follow from those of Thaddeus  so we are essentially sketching his proof.\\
We start by  introducing some notation related to the
main tool  which is a series of 
flips, depicted in the diagram (\ref{eq:z69ppr}) below, that connects the space
$\mathds{P}H^{1}(X,\xi^{-1})$
and
the moduli space $SU_{X}(2,\xi)$.

\begin{equation}\label{eq:z69ppr}
\begin{picture}(300,100)
\put(-30,20){$\mathds{P}H^{1}(X,\xi^{-1})=$}
\put(50,20){$M_{0}$}
\put(57,50){\vector(0,-1){15}}
\put(50,55){$M_{1}$}
\put(72,80){$\widetilde{M_{2}}$}
\put(75,75){\vector(-1,-1){10}}
\put(83,75){\vector(1,-1){10}}
\put(95,55){$M_{2}$}
\put(120,75){\vector(-1,-1){10}}
\put(117,80){$\widetilde{M_{3}}$}
\put(128,75){\vector(1,-1){10}}
\put(135,55){$M_{3}$}
\put(160,75){\vector(-1,-1){10}}
\put(160,80){$\widetilde{M_{4}}$} 
\put(178,75){\vector(1,-1){10}}
\put(194,55){$\dots$}
\put(222,75){\vector(-1,-1){10}}
\put(220,80){$\widetilde{M_{w}}$}
\put(230,75){\vector(1,-1){10}}
\put(237,55){$M_{w}$}
\put(243,50){\vector(0,-1){15}}
\put(234,20){$N=SU_{X}(2,\xi)$}
\end{picture}
\end{equation}
\\
The spaces $M_{i}$ are constructed as moduli spaces of  $\sigma_{i}$-stable pairs $(\mathcal{E},\phi)$ where 
$\mathcal{E}$ is a rank 2 vector bundle on the curve $X$ such that $\wedge^{2}\mathcal{E}=\xi$,  $\phi\in H^{0}(X,\mathcal{E})\backslash\{0\}$ and $\sigma_{i}\in (max(0,d/2-i-1),d/2-i)$  (for the definition of the stability condition see (1.1) in \cite{MR1273268}). 
The space $M_{1}$ is $\widetilde{\mathds{P}}^{N}_{X}$ the blow-up of 
$\mathds{P}H^{1}(X,\xi^{-1})$ along $X$ embedded via $\mid K_{X}\xi\mid$. When $i=w=[(d-1)/2]$ the map $M_{w}\mapsto N$ has fibre  $\mathds{P}H^{0}(X,\mathcal{E})$ over a stable bundle $\mathcal{E}\in N$ and it is surjective if $d>2g-2$.
All the spaces $M_{i}$  turn out to be smooth, integral, rational projective varieties of dimension $d+g-1$ and for $i>0$ there is a birational map $M_{i}\leftrightarrow M_{1}$ which is an isomorphism except on closed sets of codimension at least $2$. 
 In fact there are embeddings $\mathds{P}W_{i}^{+}\hookrightarrow M_{i}$ and 
 $\mathds{P}W_{i}^{-}\hookrightarrow M_{i-1}$, for $0<i\leq(d-1)/2$, whose images correspond to the pairs represented in $M_i$ but  not in $M_{i-1}$ and pairs represented in $M_{i-1}$ but not in $M_{i}$ respectively and there is an isomorphism
   $ M_{i}\backslash \mathds{P}W_{i}^{+}\cong M_{i-1}\backslash \mathds{P}W_{i}^{-}$.
 As for the spaces $\widetilde{M_{i}}$, for $i>1$, the arrows between $M_{i-1}$, $\widetilde{M_{i}}$ and   $M_{i}$ make up each of  the announced flips, specifically,
   $\widetilde{M_{i}}$ is the blow-up of $M_{i-1}$ along 
 $\mathds{P}W_{i}^{-}$ and it is also the the blow-up of $M_{i}$ along 
 $\mathds{P}W_{i}^{+}$,  furthermore $M_{i}$ can be obtained by  blowing-up   $M_{i-1}$ along $\mathds{P}W_{i}^{-}$ and then blowing down the same exceptional divisor to $\mathds{P}W_{i}^{+}\subset M_{i}$.
 Let $E_{i}\subset\widetilde{M_{i}}$ denote the exceptional divisor for $i=2,\cdots,w$.   If $i=1$ one can also consider 
 $\widetilde{M_{1}}\cong M_{1}$ the blow-up of $M_{1}$ along $\mathds{P}W_{1}^{+}$  with exceptional divisor $E_{1}\cong \mathds{P}W_{1}^{+}\subset M_{1}$.\\%
 For $0<i\leq w$ the embeddings $\mathds{P}W_{i}^{+}\hookrightarrow M_{i}$ and 
 $\mathds{P}W_{i}^{-}\hookrightarrow M_{i-1}$ induce  respectively the following 
 exact sequences (see (3.9) and (3.11) in \cite{MR1273268})
  \begin{equation}\label{eq:exsqtnMi}
  \begin{split}
  \mathsmall{0\rightarrow T_{\mathds{P}W_{i}^{+}}\rightarrow T_{M_{i}}\mid_{ \mathds{P}W_{i}^{+}}\rightarrow W_{i}^{-}(-1)\rightarrow 0 }~\\
  \mbox{and}\hspace{8cm}~\\
  \mathsmall{0\rightarrow T_{\mathds{P}W_{i}^{-}}\rightarrow T_{M_{i-1}}\mid_{ \mathds{P}W_{i}^{-}}\rightarrow W_{i}^{+}(-1)\rightarrow 0,}\\
\end{split}
 \end{equation}
 we recall that according to our notation here $W_{i}^{-}(-1)$  is the pull-back of $W_{i}^{-}$ to $\mathds{P}W_{i}^{+}$ under the projection to $S^{i}X$  tensored with the tautological bundle $\mathcal{O}_{\mathds{P}W_{i}^{+}}(-1)$. From the exact sequences (\ref{eq:exsqtnMi}) one notice that
 $$\mathds{P}W_{i}^{+}(-1)\cong E_{i}\cong\mathds{P}W_{i}^{-}(-1) $$ and in fact 
 $$E_{i}=\mathds{P}W_{i}^{-}\times_{S^{i}X} \mathds{P}W_{i}^{+}$$ is the fibred product of the projections  from $\mathds{P}W_{i}^{-}$ and $\mathds{P}W_{i}^{+}$ to $S^{i}X$ and we represent it in the diagram (\ref{eq:prfthad}) for future reference in the proof Lemma \ref{6.7thadppr} below.\\
 {\footnotesize
\begin{equation}\label{eq:prfthad}
\begin{picture}(390,145) 
\put(110,128){$E_{i}$}
\put(100,85){$\mathsf{p}_{1}$}
\put(116,125){\vector(0,-1){70}}
\put(130,130){\vector(1,0){95}}
\put(170,135){$\mathsf{p}_{2}$}
\put(230,128){$\mathds{P}W_{i}^{+}$}
\put(240,125){\vector(0,-1){70}}

\put(245,85){$\pi_{2}$}
\put(230,45){$S^{i}X,$}
\put(105,45){$\mathds{P}W_{i}^{-}$}
\put(130,45){\vector(1,0){95}}
\put(170,50){$\pi_{1}$}
\end{picture}
\end{equation}
}
One denotes by $\mathcal{O}_{E_{i}}(j,k)$ the line bundle on $E_{i}$ given by  
\begin{equation}\label{eq:OEijk}
\mathcal{O}_{E_{i}}(j,k):=\mathcal{O}_{\mathds{P}W_{i}^{-}}(j)\boxtimes \mathcal{O}_{\mathds{P}W_{i}^{+}}(k).
\end{equation}
The twisting sheaf of $E_{i}$ happens to be
$$\mathcal{O}_{E_{i}}(1)=\mathcal{O}_{E_{i}}(1,1)$$ 
therefore
$$\mathcal{O}_{E_{i}}(E_{i})=\mathcal{O}_{E_{i}}(-1,-1).$$

~ 

 Now we consider the line bundles $\mathcal{O}_{1}(m,n)=\mathcal{O}_{M_{1}}((m+n)H-nE)$ over $M_{1}$. 
Since $M_{1}$ and $M_{i}$ are isomorphic outside  closed sets of codimension at least $2$
there is (as a consequence of Hartogs' theorem), for $i>0 $,  an isomorphism
\begin{equation*}\label{eq:82ppr}
Pic \ M_{1}\cong Pic \ M_{i}. 
\end{equation*}
Let $\mathcal{O}_{i}(m,n)$ be the image of $\mathcal{O}_{1}(m,n)$ under this isomorphism. 
One also has (another consequence of Hartogs' theorem) that
\begin{equation}\label{eq:99ppr}
H^{0}(M_{1},\mathcal{O}_{1}(m,n))=H^{0}(M_{i},\mathcal{O}_{i}(m,n)).
\end{equation}
Recall that we denoted in the introduction $V_{m,n}= H^{0}(M_{1},\mathcal{O}_{1}(m,n)$.
Abusing notation one also writes  $\mathcal{O}_{i}(m,n)$ and $\mathcal{O}_{i-1}(m,n)$ for the pull-backs  to the space $\widetilde{M_{i}}$ of the line bundles  $\mathcal{O}_{i}(m,n)$ and $\mathcal{O}_{i-1}(m,n)$ defined above.
In the case $i=0$  we have $M_{0}=\mathds{P}H^{1}(X,\xi^{-1})$ and one  defines  $\mathcal{O}_{0}(m,n)=\mathcal{O}_{\mathds{P}H^{1}(X,\xi^{-1})}(m+n)$.
 On 
 $\widetilde{M_{i}}$ one has (see (5.6) in \cite{MR1273268})
 \begin{equation}\label{eq:thad5.6}
 \mathcal{O}_{i}(m,n)=\mathcal{O}_{i-1}(m,n)(((i-1)m-n)E_{i}).
 \end{equation}
 
When we assume that a group $G$ is acting on the curve $X$ one has by functoriality an action on all the moduli spaces $M_{i}$ and $N$. This action is compatible with the isomorphisms
$ M_{i}\backslash \mathds{P}W_{i}^{+}\cong M_{i-1}\backslash \mathds{P}W_{i}^{-}$ because pairs representing points on the left-hand side also represent points on the right-hand side.
In the introduction we endowed the line bundles $\mathcal{O}_{1}(m,n)$ with a G-linearization and because of the isomorphism 
$M_{1}\leftrightarrow M_{i}$
 outside closed subsets of codimension at least $2$ one has that
this linearization also induces a linearization on the other line bundles  $\mathcal{O}_{i}(m,n)$. In this way, the action of $G$ induced on 
$H^{0}(M_{1},\mathcal{O}_{1}(m,n))$ and on $H^{0}(M_{i},\mathcal{O}_{i}(m,n))$ through these linearizations is the same. The action of $G$ on the spaces $M_{i}$ also lifts to an action on the spaces $\widetilde{M_{i}}$ and one also  endows the line bundles $\mathcal{O}_{i}(m,n),\mathcal{O}_{i-1}(m,n)$ in $Pic(\widetilde{M_{i}})$ with a linearization induced by that of $\mathcal{O}_{i}(m,n)$ in  $Pic(M_{i})$ and $\mathcal{O}_{i-1}(m,n)$ in  $Pic(M_{i-1})$.  Now we can state the following lemmas.
\begin{l1}\label{6.2thadppr}If $m,n\geq 0$ and $m(d-2)-2n>-d+2g-2$ there exist an  integer $b\leq w$ such that  for any $h\in G$
\begin{equation*}\label{eq:2.6vppr}
Trace(h_{\mid} V_{m,n})=L(h,M_{b},\mathcal{O}_{b}(m,n)),
\end{equation*}
in fact, $b=[\frac{n+d+g-4}{m+3}]+1.$ 
\end{l1}
\begin{proof}
Thaddeus showed that for $b=[\frac{n+d+g-4}{m+3}]+1$
the  cohomology groups $H^{i}(M_{b},\mathcal{O}_{b}(m,n))$ vanish for $i>0$ (see proof of (6.2) in \cite{MR1273268}). 
Then for any $h\in G$
\begin{multline*}
L(h,M_{b},\mathcal{O}_{b}(m,n))=\sum_{i\geq 0}(-1)^{i}Trace \ h_{|}H^{i}(M_{b},\mathcal{O}_{b}(m,n))\\
\hspace{3cm}~
=Trace \ h_{|}H^{0}(M_{b},\mathcal{O}_{b}(m,n)).\hspace{3.5cm}~
\end{multline*}
Using (\ref{eq:99ppr}) we have
\[ 
Trace \ h_{|}H^{0}(M_{b},\mathcal{O}_{b}(m,n))=Trace \ h_{|}H^{0}(M_{1},\mathcal{O}_{1}(m,n)) 
\]
and since $V_{m,n}=H^{0}(M_{1},\mathcal{O}_{1}(m,n))$ we are done.
\end{proof}
\noindent
For the remainder of this section we shall assume that
$m,n\geq 0, \ m(d-2)-2n>-d+2g-2$ and $b=[\frac{n+d+g-4}{m+3}]+1.$
\begin{l1}\label{6.6thadppr}
$N_{0}(h)=L(h,M_{0},\mathcal{O}_{0}(m,n)).$
\end{l1}
\begin{proof} This follows from the  definitions of  $N_{0}(h)$ and $\mathcal{O}_{0}(m,n).$ Namely, we have  $M_{0}=\mathds{P}H^{1}(X,\xi^{-1})$ and   $\mathcal{O}_{0}(m,n)=\mathcal{O}_{\mathds{P}H^{1}(X,\xi^{-1})}(m+n)$ (see definition before (\ref{eq:thad5.6})). Notice that when $i=0$,  we have $q_{i}-i=m+n\geq 0$ (see definition before (\ref{eq:z57})). So we have the following
$$H^{0}(M_{0},\mathcal{O}_{0}(m,n))=S^{m+n}H^{0}(\mathds{P}H^{1}(X,\xi^{-1}),\mathcal{O}_{\mathds{P}H^{1}(X,\xi^{-1})}(1))$$
$$=S^{m+n}H^{0}(X,K_{X}\xi)$$
and 
$H^{i}(M_{0},\mathcal{O}_{0}(m,n))=0$ for $i>0$. Then 
$$L(h,M_{0},\mathcal{O}_{0}(m,n))=Trace(h_{\mid_{}}S^{m+n}H^{0}(X,K_{X}\xi))$$
and by the definition of $B_{0,m,n}$ in Section \ref{stf} the last is
$$=Trace(h_{\mid_{}}B_{0,m,n})$$
which in turn is the definition of $N_{0}(h)$
in  (\ref{eq:3c2p2i0}) .
\end{proof}
\noindent We  can also restate  (6.7) of \cite{MR1273268} as the following 
\begin{l1}\label{6.7thadppr} For $0<i\leq b$ we have for any $h\in G$
\begin{equation*}
L(h,M_{i},\mathcal{O}_{i}(m,n))-L(h,M_{i-1},\mathcal{O}_{i-1}(m,n))=(-1)^{i}N_{i}(h)
\end{equation*}
\end{l1}
\begin{proof}
When $h$ is the identity Thaddeus first shows that 
$\chi(M_{i},\mathcal{O}_{i}(m,n))=\chi(\widetilde{M_{i}},\mathcal{O}_{i}(m,n))$. But this follows because   the  higher direct images of $\mathcal{O}_{i}(m,n)$ under the map $\pi:\widetilde{M_{i}}\mapsto M_{i}$ vanish so the Leray spectral sequence  implies 
$H^j(\widetilde{M_{i}},\mathcal{O}_{i}(m,n))=H^j(M_{i},R^0_{\pi*}\mathcal{O}_{i}(m,n))$
 and then since $R^0_{\pi*}\mathcal{O}_{\widetilde{M_{i}}}=\mathcal{O}_{M_{i}}$ the projection formula implies
 $H^j(\widetilde{M_{i}},\mathcal{O}_{i}(m,n))=H^j(M_{i},\mathcal{O}_{i}(m,n))$. Similarly $H^j(\widetilde{M_{i}},\mathcal{O}_{i-1}(m,n))=H^j(M_{i-1},\mathcal{O}_{i-1}(m,n))$. Then we can write for any $h\in G$
 $$L(h,M_{i},\mathcal{O}_{i}(m,n))=L(h,\widetilde{M_{i}},\mathcal{O}_{i}(m,n))$$
 and
 $$L(h,M_{i-1},\mathcal{O}_{i-1}(m,n))=L(h,\widetilde{M_{i}},\mathcal{O}_{i-1}(m,n))$$
 and one can work on $\widetilde{M_{i}}$.\\
 Recall that $q_{i}=n-(i-1)m$. Next one considers 2 cases.
\\1) When $q_{i}\leq 0$,  by definition of $B_{i,m,n}$ (see equation (\ref{eq:z58})),  one has $N_{i}(h)=0$.\\
If $q_{i}= 0$ then we are done because  $\mathcal{O}_{i}(m,n)=\mathcal{O}_{i-1}(m,n)$ on $\widetilde{M_{i}}$ by  (\ref{eq:thad5.6}). So we assume now that $q_{i}< 0$.\\
  Since $E_{i}$ is a $G$-invariant divisor on $\widetilde{M_{i}}$ the exact sequence induced by the embedding $E_{i}\hookrightarrow \widetilde{M_{i}}$
 \begin{equation}
  \mbox{\small$0\rightarrow\mathcal{O}_{\widetilde{M_{i}}}(-E_{i})\rightarrow\mathcal{O}_{\widetilde{M_{i}}}\rightarrow\mathcal{O}_{\widetilde{M_{i}}}\otimes\mathcal{O}_{E_{i}}\rightarrow 0 $}
 \end{equation}
 is a $G$-equivariant exact sequence of linearized sheaves so 
     for each $ j$ we have an equivariant exact sequence
 \begin{equation}\label{eq:es1pr6.1}
  \mbox{\footnotesize$0\rightarrow\mathcal{O}_{i-1}(m,n)((j-1)E_{i})\rightarrow\mathcal{O}_{i-1}(m,n)(jE_{i})\rightarrow\mathcal{O}_{i-1}(m,n)\otimes\mathcal{O}_{E_{i}}(jE_{i})\rightarrow 0. $}
 \end{equation}
 Thaddeus identify $\mathcal{O}_{i-1}(m,n)\otimes\mathcal{O}_{E_{i}}(jE_{i})$ with the sheaf(or rather its pushforward to $\widetilde{M_{i}}$) 
 $L_{i}^{m}(-q_{i}-j,-j)$.  The last corresponds to a the tensor product $F_{1}\otimes F_{2}$ of sheaves on $E_{i}$  where $F_{1}$ is the pull-back of the sheaf $L_{i}^{m}$ ( defined in (\ref{eq:z56})) under the projection
 $\pi:\mathftnt{E_{i}\twoheadrightarrow S^{i}X}$
  (here $\pi$ is $\pi_{1}\circ\mathsf{p}_{1}= \pi_{2}\circ\mathsf{p}_{2}$ in  the diagram (\ref{eq:prfthad}))
  and $F_{2}=\mathcal{O}_{E_{i}}(-q_{i}-j,-j)$ (see equation (\ref{eq:OEijk})).\\
    The  exact sequence (\ref{eq:es1pr6.1}) induces  a $G$-equivariant long exact sequence of cohomology groups from which we can write ( using the above identification):
\footnotesize
\begin{multline}\label{eq:ps1prtd6.7}
\!\!\!\!L(h,\widetilde{M_{i}},\mathcal{O}_{{}_{i-1}}\!(m,n)(jE_{i}))-L(h,\widetilde{M_{i}},\mathcal{O}_{{}_{i-1}}\!(m,n)(\!(j-1)E_{1}\!))=
\\
=
L(h,E_{i},L_{i}^{m}(\!-q_{i}\!-\!j,\!-\!j)).\mbox{{ }  }\!\!\!\!
\end{multline}
\normalsize
 Summing (\ref{eq:ps1prtd6.7}) over $j$, with $0<j\leq-q_{i}$, and using (\ref{eq:thad5.6}) one arrives to
 \begin{equation}\label{eq:s2prtd6.7}
 \mbox{\footnotesize$L(h,\widetilde{M_{i}},\mathcal{O}_{i}(m,n))-L(h,\widetilde{M_{i}},\mathcal{O}_{i-1}(m,n))=\sum_{j=1}^{-q_{i}}L(h,E_{i},L_{i}^{m}(-q_{i}-j,-j)).$}
 \end{equation}
 When $h$ is the identity Thaddeus proves the vanishing of the right-hand side by showing that all  the cohomology groups 
$$H^{s}(E_{i},L_{i}^{m}(-q_{i}-j,-j))$$
vanish (all the direct images of  $L_{i}^{m}(-q_{i}-j,-j)$ under the projection $\mathsf{p}_{1}:E_{i}\rightarrow \mathds{P}W_{i}^{-}$ vanish because 
$\mbox{\footnotesize$R^{t}_{\pi*}L_{i}^{m}(-q_{i}-j,-j)=L_{i}^{m}(-q_{i}-j)\otimes R^{t}_{\pi*}\mathcal{O}_{E_{i}}(0,-j)$}$ 
\\
and  
$\mbox{\footnotesize $R^{t}_{\pi*}\mathcal{O}_{E_{i}}(0,-j)=0$}$ for all $t$ because $0<j<d+g-1-2i=\mbox{rank$(W_{i}^{+})$}$ implies that a fibre
$\mathftnt{H^t( \mathds{P}^{d+g-2-2i},\mathcal{O}(-j))=0}$).
 Therefore the right-hand side of (\ref{eq:s2prtd6.7}) vanish for any $h\in G$.\\
2) Now we assume $q_{i}>0.$ Summing (\ref{eq:ps1prtd6.7}) over $j$, with $-q_{i}< j\leq 0$, and using (\ref{eq:thad5.6}) one obtains
\begin{equation}\label{eq:s3prtd6.7}
\begin{split}
 \mbox{\footnotesize$L(h,\widetilde{M_{i}},\mathcal{O}_{i-1}(m,n))-L(h,\widetilde{M_{i}},\mathcal{O}_{i}(m,n))$}&=\sum_{j=-q_{i}+1}^{0}\!\!\mbox{\footnotesize$L(h,E_{i},L_{i}^{m}(-q_{i}-j,-j))$}\\
 &= \sum^{q_{i}-1}_{j=0}\mbox{\footnotesize$L(h,E_{i},L_{i}^{m}(-q_{i}+j,+j)).$}\\
 \end{split}
 \end{equation}
Now, the only non-zero direct image of  $\mathftnt{L_{i}^{m}(-q_{i}+j,+j)}$  under the projection $\pi:\mathftnt{E_{i}\twoheadrightarrow S^{i}X}$
 is the (i-1)-th (see (\ref{eq:drimgi-1}) below):\\
  given $D\in S^{i}X$ a fibre of $\pi$ is of the form ${E_{i}}_{D}=\mathds{P}^{i-1}\times \mathds{P}^{d+g-2-2i}$,\\
  since $-q_{i}+j<0$ one has $H^{s}(\mathds{P}^{i-1},\mathcal{O}(-q_{i}+j))=0$ for $s\not=i-1$ and\\
  since  $j\geq0$, $H^{s}(\mathds{P}^{d+g-1-2i},\mathcal{O}(+j))=0$ for $s\not=0$.\\
  So $H^{i-1}(\mathds{P}^{i-1},\mathcal{O}(-q_{i}+j))\otimes H^{0}(\mathds{P}^{d+g-1-2i},\mathcal{O}(+j))$ is the only non-zero  K\"{u}n\-neth component of 
  $H^{i-1}({E_{i}}_{D},\mathcal{O}_{E_{i}}(-q_{i}+j,+j)\mid_{{E_{i}}_{D}})$. \\
  To compute a direct image ${R^{i-1}_{\pi*}}_{}F$  one uses the fact that the Leray Spectral sequence $E_{2}^{t,s}={R^{t}_{\pi_{1}*}}_{}{R^{s}_{\mathsf{p}_{1}*}}_{}F\Rightarrow$
  ${R^{t+s}_{\pi*}}_{}F$. It will be enough to work with $F=$ $\mathcal{O}_{E_{i}}(-q_{i}+j,+j)$ because $R^{i-1}_{\pi*}(L_{i}^{m}(-q_{i}+j,+j)) = L_{i}^{m}\otimes R^{i-1}_{\pi*}\mathcal{O}_{E_{i}}(-q_{i}+j,+j)$. 
  We shall see that the only non-zero term of the Leray spectral sequence is $E_{2}^{i-1,0}$.\\
  Notice that for all sheaves ${\mathcal{O}}_{\mathds{P}W_{i}^{-}}(l)$, ${\mathcal{O}}_{\mathds{P}W_{i}^{+}}(l)$ the base change morphisms induced by the diagram (\ref{eq:prfthad}) are isomorphisms, that is,\\
  ${R^{t}_{\mathsf{p}_{1}*}}_{}\mathsf{p}_{2}^{*}{\mathcal{O}}_{\mathds{P}W_{i}^{+}}(l)\cong \pi_{1}^{*}{R^{t}_{\pi_{2}*}}_{}{\mathcal{O}}_{\mathds{P}W_{i}^{+}}(l)$
  and
  ${R^{t}_{\mathsf{p}_{2}*}}_{}\mathsf{p}_{1}^{*}{\mathcal{O}}_{\mathds{P}W_{i}^{-}}(l)\cong \pi_{2}^{*}{R^{t}_{\pi_{1}*}}_{}{\mathcal{O}}_{\mathds{P}W_{i}^{-}}(l)$ for all $t\geq 0$.\\
  Now, by the projection formula one has\\
  ${R^{t}_{\pi_{1}*}}_{}{R^{s}_{\mathsf{p}_{1}*}}_{}\mathcal{O}_{E_{i}}(-q_{i}+j,+j)=$ 
  ${R^{t}_{\pi_{1}*}}_{}\{ \mathcal{O}_{\mathds{P}W_{i}^{-}}(-q_{i}+j)\otimes {R^{s}_{\mathsf{p}_{1}*}}_{}\mathsf{p}_{2}^{*}{\mathcal{O}}_{\mathds{P}W_{i}^{+}}(+j)\}=$\\
  (by the base change isomorphism) 
 ${R^{t}_{\pi_{1}*}}_{}\{ \mathcal{O}_{\mathds{P}W_{i}^{-}}(-q_{i}+j)\otimes \pi_{1}^{*}{R^{s}_{\pi_{2}*}}_{}{\mathcal{O}}_{\mathds{P}W_{i}^{+}}(+j)\}=$ \\
 ${R^{t}_{\pi_{1}*}}_{}\{ \mathcal{O}_{\mathds{P}W_{i}^{-}}(-q_{i}+j)\}\otimes {R^{s}_{\pi_{2}*}}_{}\{{\mathcal{O}}_{\mathds{P}W_{i}^{+}}(+j)\}.$ \\
 Since  $j\geq0$, for $s>0$ one has ${R^{s}_{\pi_{2}*}}_{}\{{\mathcal{O}}_{\mathds{P}W_{i}^{+}}(+j)\}=0$
 and for $s=0$ ${R^{s}_{\pi_{2}*}}_{}\{{\mathcal{O}}_{\mathds{P}W_{i}^{+}}(+j)\}=$
 $ S^{j}(W_{i}^{+})^{\vee}.$\\
 Since   $-q_{i}+j<0$ we have that ${R^{t}_{\pi_{1}*}}_{}\{ \mathcal{O}_{\mathds{P}W_{i}^{-}}(-q_{i}+j)\}=0$ 
 unless $t= \mbox{rank}(W_{i}^{-})-1 = i-1$ in which case \\
 ${R^{i-1}_{\pi_{1}*}}_{}\{ \mathcal{O}_{\mathds{P}W_{i}^{-}}(-q_{i}+j)\}=\wedge^{i}W_{i}^{-}\otimes S^{q_{i}-j-i}(W_{i}^{-}).$
 
\noindent It follows that
\begin{equation}\label{eq:drimgi-1}
\mathftnt{R^{i-1}_{\pi*}(L_{i}^{m}(-q_{i}+j,+j))=L_{i}^{m}\otimes\wedge^{i}W_{i}^{-}\otimes S^{q_{i}-j-i}(W_{i}^{-})\otimes S^{j}(W_{i}^{+})^{\vee}}.
\end{equation}
 By the Leray spectral sequence(the usual one) one has
\begin{equation}
\mathftnt{H^{s}(E_{i},L_{i}^{m}(-q_{i}+j,+j))=H^{s-(i-1)}(S^{i}X,R^{i-1}_{\pi*}(L_{i}^{m}(-q_{i}+j,+j)))}
\end{equation}
and the rest is a rather straightforward verification. One can write
\begin{equation}\label{eq:idbthad2.6}
\begin{split}
\mathftnt{L(h,E_{i},L_{i}^{m}(-q_{i}+j,+j))}= {\:\:\:\:\:\:\:\:\:\:\:\:\:\:\:\:\:\:\:\:\:\:\:\:\:\:\:\:\:\:\:\:\:\:\:\:\:\:\:\:\:\:\:\:\:\:\:\:\:\:\:\:\:\:\:\:\:\:\:\:\:\:\:\:\:\:\:\:\:\:\:\:\:\:\:\:\:\:\:\:}\\
{\:\:\:\:\:\:\:\:\:\:\:\:\:\:\:\:\:\:\:\:\:\:\:\:\:\:\:\:\:\:\:\:}\sum_{s}\mathftnt{(-1)^{s}trace\left( h_{\mid}H^{s-(i-1)}(S^{i}X,R^{i-1}_{\pi*}(L_{i}^{m}(-q_{i}+j,+j)))\right)}\\
{\:\:\:\:\:\:\:\:\:\:\:\:\:\:\:\:\:\:\:\:\:\:\:\:\:\:\:\:\:\:\:\:}=(-1)^{i-1}\mathftnt{L\left( h,S^{i}X,R^{i-1}_{\pi*}(L_{i}^{m}(-q_{i}+j,+j))\right)}.\\
\end{split}
\end{equation}
Multiplying  (\ref{eq:s3prtd6.7}) by $-1$  and using (\ref{eq:idbthad2.6}) we have
\begin{multline*}
\mathftnt{L(h,\widetilde{M_{i}},\mathcal{O}_{i}(m,n))-L(h,\widetilde{M_{i}},\mathcal{O}_{i-1}(m,n))}=\\
\hspace{2.cm}~-\left[\mathftnt{L(h,\widetilde{M_{i}},\mathcal{O}_{i-1}(m,n))-L(h,\widetilde{M_{i}},\mathcal{O}_{i}(m,n))}\right]=\hspace{2cm}~\\
\hspace{2.cm}~=(-1)^{i}\sum_{j=0}^{q_{i}-1}\mathftnt{L\left( h,S^{i}X,R^{i-1}_{\pi*}(L_{i}^{m}(-q_{i}+j,+j))\right)}\\
\end{multline*}
from (\ref{eq:drimgi-1}) this is
\begin{multline*}
=(-1)^{i}\sum_{j=0}^{q_{i}-1}\mathftnt{L\left( h,S^{i}X,L_{i}^{m}\otimes\wedge^{i}W_{i}^{-}\otimes S^{q_{i}-j-i}(W_{i}^{-})\otimes S^{j}(W_{i}^{+})^{\vee}\right)}\\
\end{multline*}
 and if $j>q_{i}-i$ then $S^{q_{i}-j-i}(W_{i}^{-})=0$ so in the last equation the sum becomes 
a sum that runs only from $j=0$ to $j=q_{i}-i$, that is
\begin{multline*}
=(-1)^{i}\sum_{j=0}^{q_{i}-i}\mathftnt{L\left( h,S^{i}X,L_{i}^{m}\otimes\wedge^{i}W_{i}^{-}\otimes S^{q_{i}-j-i}(W_{i}^{-})\otimes S^{j}(W_{i}^{+})^{\vee}\right)}{}\\
=\mathftnt{(-1)^{i}L}\left({\mathftnt{ h,S^{i}X,L_{i}^{m}\otimes\wedge^{i}W_{i}^{-}\otimes}\left\lbrace\bigoplus_{j=0}^{q_{i}-i} \mathftnt{ S^{q_{i}-j-i}(W_{i}^{-})\otimes S^{j}(W_{i}^{+})^{\vee}}\right\rbrace}\right)\hspace{1cm}~\\
(\mbox{using } \bigoplus_{j=0}^{q_{i}-i} {\mathftnt{ S^{q_{i}-j-i}(W_{i}^{-})\otimes S^{j}(W_{i}^{+})^{\vee}}}=S^{q_{i}-i}\left\lbrace {\mathftnt{ W_{i}^{-}\oplus (W_{i}^{+})^{\vee}}}\right\rbrace)\\
=\mathftnt{(-1)^{i}L}\left({\mathftnt{ h,S^{i}X,L_{i}^{m}\otimes\wedge^{i}W_{i}^{-}\otimes S^{q_{i}-i}}\left\lbrace \mathftnt{ W_{i}^{-}\oplus (W_{i}^{+})^{\vee}}\right\rbrace}\right)\hspace{3cm}~\\
=\mathftnt{(-1)^{i}L}\left({\mathftnt{ h,S^{i}X,L_{i}^{m}\otimes\wedge^{i}W_{i}^{-}\otimes S^{q_{i}-i} U_{i}}}\right){\:\:\:\:\:\:\:\:\:\:\:\:\:\:\:\:\:}\\
=\mathftnt{(-1)^{i}L}\left({\mathftnt{ h,S^{i}X,B_{i,m,n}}}\right){\:\:\:\:\:\:\:\:\:\:\:\:\:\:\:\:\:}\\
=\mathftnt{(-1)^{i}N_{i}(h)}.{\:\:\:\:\:\:\:\:\:\:\:\:\:\:\:\:\:}\\
\end{multline*}
\end{proof}
\begin{l1}\label{6.8thadppr}
For $i>b$, $N_{i}(h)=0$.
\end{l1}
\begin{proof}
This is the same proof of (6.8) in \cite{MR1273268} since one shows that $i>b$ implies $q_{i}-i<0$ so $B_{i,m,n}=0$.
\end{proof}
\noindent Finally, from Lemmas \ref{6.2thadppr}--\ref{6.8thadppr} one has\\
 $\sum_{i=1}^{w}(-1)^{i}N_{i}(h)=L(h,M_{b},\mathcal{O}_{b}(m,n))-L(h,M_{0},\mathcal{O}_{0}(m,n))$\\
$~ \hspace{3cm} =Trace( h_{\mid_{ V_{m,n}}} )-N_{0}(h).$\\
Therefore $Trace( h_{\mid_{ V_{m,n}}} )=\sum_{i=0}^{w}(-1)^{i}N_{i}(h)$.
 \section{The Chern classes}\label{sec:TchC}
Let $h$ be an automorphism of the curve $X$ and assume that  $h$ has order $p\not= 1$.
We shall explain below that the $k$-dimensional components of fixed points of $h$ in the symmetric product $S^{i}X$ are parametrized by certain kind of $h$-invariant divisors so we represent such a component by $Z_{D}$, where $D$ is the corresponding invariant 
divisor.
Let $\iota_{D}$ denote the inclusion  $Z_{D}\subset S^{i}X$. Consider the decompositions into  eigenbundles 
\begin{equation}\label{eq:z20}
\iota_{D}^{*}W^{-}_{i}=\bigoplus_{j=1}^{p} \iota_{D}^{*}W^{-}_{i}(\nu	^{j})
\end{equation}
and
\begin{equation}\label{eq:z21}
\iota_{D}^{*}W^{+}_{i}=\bigoplus_{j=1}^{p} \iota_{D}^{*}W^{+}_{i}(\nu	^{j}).
\end{equation}

For the proof of Theorem \ref{t:chhBimn} in Section \ref{TgCcfBimn} we  need to know the Chern classes of all these eigenbundles and  before we compute them we recall from \cite{MR2125533} Section 3,  that  a $k$-dimensional component of fixed points of $h$ in $S^{i}X$ is isomorphic to the symmetric product $S^{k}Y$, here $Y$ is the quotient curve $X/\<h\>$. These components are parametrized by a set of certain kind of  $h$-invariant divisors $A_{k}$ of degree $d_{k}=i-pk$, hence the notation $Z_{D}$. More precisely, define $A_{k}$ as the set of divisors $D\in {(S^{d_{k}}X)}^{h}$ satisfying the following property: if $x\in X$ is a point in the support of $D$ then $D-\sum_{j=0}^{p-1}h^{j}x$ is not an effective divisor nor the zero divisor.
For each $D\in A_{k}$ there is an embedding
\begin{equation}\label{eq:degreecomposition}
\iota_{D}:S^kY\stackrel{\iota}{\hookrightarrow} 
S^{pk}X\stackrel{\mathcal{A}_{D}}{\hookrightarrow} S^{pk+d_{k}}X 
\end{equation}
 where $\iota$ sends $P\in S^kY$ to the divisor $f^*P\in 
S^{pk}X$ ( $f:X\rightarrow Y=X/\<h\>$ is the quotient map) and 
$\mathcal{A}_{D}$ 
sends 
$P\in S^{pk}X$ to $P+D\in S^{pk+d_{k}}X$. \\
Then the Chern classes of our eigenbundles can be expressed in terms of the cohomology classes $\theta,x$ and $\sigma_{i}\in H^2(S^kY,\mathbb{Z})$  (see \cite{MR151460} for details on cohomology of symmetric products), where $x$ represents the class of a divisor $q+ S^{k-1}Y\subset S^kY$ in $H^2(S^kY,\mathbb{Z})$ and  $\theta $ represents the class of the pull back of the theta divisor class 
$\Theta\in H^2(J_{Y},\mathbb{Z})$ of the Jacobian $J_{Y}$ of the curve $Y$ under the Abel Jacobi map. We recall some relations of these  cohomology classes:
\begin{equation}\label{eq:z21theta}
\theta=\sum_{i=1}^{g_{Y}}\sigma_{i}, \mbox{ }\sigma_{i}\sigma_{j}= \sigma_{j}\sigma_{i}, \mbox{ and } \sigma_{i}^{2}=0.
\end{equation}
If $0\leq a\leq g_{Y}$ and $0\leq d$, then 
\begin{equation}
\sigma_{i_{1}}\sigma_{i_{2}}\dots \sigma_{i_{a}}x^{d}=x^{a+d} , \mbox{ for distinct $i_{1}, i_{2},\dots ,i_{a}$. Also}
\end{equation}

\begin{equation}
\theta^{a}x^{d}=a!\binom{g_{Y}}{a}x^{a+d} \mbox{ and}
\end{equation}
\begin{equation}\label{eq:z21exptheta}
e^{z\theta}=\prod_{i=1}^{g_{Y}}(1+z\sigma_{i}).
\end{equation}
The  Todd  class $Td(Z_{D})$ of a k-dimensional component is given by 
\begin{equation}\label{eq:97}
Td(Z_{D})=Td(S^{k}Y)=\left(\frac{x}{1-e^{-x}}\right)^{k-g_{Y}+1}exp\left(\frac{\theta}{e^{x}-1}-\frac{\theta}{x}\right), 
\end{equation}
 see for instance (7.3) in  \cite{MR1273268}, there $\sigma=\theta$ and $\eta=x$.

Let $D\in A_{k}$.  We will consider as well the following decomposition into eigenbundles on Y
\begin{equation}\label{eq:z22n}
f_{*}(\xi^{s}(-nD))=\bigoplus^{p-1}_{j=0}\lambda_{s,j,n},
\end{equation}
\nothing{ 
\begin{equation}\label{eq:z22}
f_{*}(\xi\otimes \mathcal{O}_{X}(-nD))=\bigoplus^{p-1}_{j=0}\mathcal{L}_{j,n}
\end{equation}    
and
\begin{equation}\label{eq:z39}
f_{*}(\xi^{-1}\otimes \mathcal{O}_{X}(-nD))=\bigoplus^{p-1}_{j=0}\mathcal{L'}_{j,n}.
\end{equation} 
}
here $\lambda_{s,j,n}:=f_{*}(\xi^{s}(-nD))(\nu^{j})$.
We have the following
\begin{t1}\label{chfw+w-} Let $Z_{D}$ be a k-dimensional component of fixed points of $h$ in $S^{i}X$. Let $m_{j,1},m_{j,2},$$m'_{j,n}$ denote the degrees of the bundles $\lambda_{1,j,1}$,$\lambda_{1,j,2}$ and $\lambda_{-1,j,n}$,respectively, in formula (\ref{eq:z22n}). Let $g_{Y}$ be the genus of the quotient curve $Y$. Then,\\
a) For the eigenbundles in (\ref{eq:z21})  their corresponding Chern  characters and  classes are given by:

\begin{equation}\label{eq:z37}
ch(\iota_{D}^{*}W^{+}_{i}(\nu^{j}))=-e^{2x}(1+m'_{j,-2}-(-2k+g_{Y}+4\theta))
\end{equation}
and
\begin{equation}\label{eq:z38}
c(\iota_{D}^{*}W^{+}_{i}(\nu^{j}))=\frac{e^{\frac{4\theta}{1+2x}}}{(1+2x)^{(1+m'_{j,-2}+2k-g_{Y})}}.
\end{equation}
b) For the eigenbundles in equation (\ref{eq:z20}) we have:
\begin{equation}
ch(\iota_{D}^{*}W^{-}_{i}(\nu^{j}))=e^{-x}(1+m_{j,1}-(k+g_{Y}+\theta))-e^{-2x}(1+m_{j,2}-(2k+g_{Y}+4\theta))
\end{equation}
and
\begin{equation}
c(\iota_{D}^{*}W^{-}_{i}(\nu^{j}))=\frac{(1-x)^{1+m_{j,1}-k-g_{Y}}}{(1-2x)^{1+m_{j,2}-2k-g_{Y}}}
e^{-\frac{\theta}{1-x}+\frac{4\theta}{1-2x}}.
\end{equation}

\end{t1}
  In the diagrams (\ref{eq:z23}),(\ref{eq:z267aart})and (\ref{eq:z268art}) below we introduce  notation for some morphisms that appear in the proof of Theorem \ref{chfw+w-} and of Lemma \ref{pLchfw+w-}.  
In the  diagram (\ref{eq:z23})  $\rho_{S^{k}Y}$ and $\pi_{S^i X}$ are the natural projections, $\iota_{D}$ is the embedding  (\ref{eq:degreecomposition}) corresponding to the component $Z_{D}$, $j'$ stands for the embedding of the universal divisor $\Delta_{i}$ of $S^{i}X$,  
\[
\LinL{\Delta}:= (Id_{X}\times\iota_{D})^{*} \Delta_{i}
\]
and $\beta'$ is the corresponding embedding into $X\times S^{k}Y$.
{\footnotesize
\begin{equation}\label{eq:z23}
\begin{picture}(390,230) 

\put(83,210){$\LinL{\Delta}$}
\put(85,180){
\begin{tikzpicture}
\pgfsetarrows{inyectivo2->}
\pgfsetlinewidth{.1ex}
\pgfpathmoveto{\pgfpointorigin}
\pgfpathlineto{\pgfpoint{0cm}{-.85cm}}
\pgfusepath{stroke}
\end{tikzpicture}}
\put(69,193){$\beta'$}
\nothing{
\put(195,201){$\phi$}
\put(98,212){\vector(1,0){198}}
}
\nothing{
\put(125,200){$\xi''$}
\put(125,198){\vector(-2,-1){24}}
} 
\put(70,173){$X\times S^{k}Y$}
\put(55,120){$\rho_{S^{k}Y}$}
\put(86,170){\vector(0,-1){119}}
\put(120,175){\vector(1,0){165}}
\put(190,180){$Id_{X}\times\iota_{D}$}
\nothing{
\put(180,40){$S^{pk}X$}

\put(205,43){
\begin{tikzpicture}
\pgfsetarrows{inyectivo2-arcs'}
\pgfsetlinewidth{.1ex}
\pgfpathmoveto{\pgfpointorigin}
\pgfpathlineto{\pgfpoint{2.8cm}{0cm}}
\pgfusepath{stroke}
\end{tikzpicture}
}

\put(240,48){$\gamma$}
} 
\put(290,173){$X\times S^{i}X$}
\put(310,172){\vector(0,-1){121}}

\put(318,120){$\pi_{S^{i}X}$}
\nothing{
\put(333,121){$X$}
\put(318,170){\vector(1,-2){20}}
\put(329,150){$P$}
} 
\nothing{
\put(255,200){$\xi_{S^{pk+d_{k}}X}$}
\put(280,197){\vector(3,-2){20}}
}
\put(305,210){$\Delta_{i}$}
\put(307,180){
\begin{tikzpicture}
\pgfsetarrows{inyectivo2->}
\pgfsetlinewidth{.1ex}
\pgfpathmoveto{\pgfpointorigin}
\pgfpathlineto{\pgfpoint{0cm}{-.85cm}}
\pgfusepath{stroke}
\end{tikzpicture}}
\put(318,195){$j'$}
\put(303,40){$S^{i}X,$}
\put(75,40){$S^{k}Y$}
\put(100,45){\vector(1,0){200}}
\put(200,50){$\iota_{D}$}
\nothing{
\put(95,43){
\begin{tikzpicture}
\pgfsetarrows{inyectivo2-arcs'}
\pgfsetlinewidth{.1ex}
\pgfpathmoveto{\pgfpointorigin}
\pgfpathlineto{\pgfpoint{2.8cm}{0cm}}
\pgfusepath{stroke}
\end{tikzpicture}}

\put(140,48){$\beta$}
}

\nothing{\put(122,70){
\begin{tikzpicture}[line width=.5pt]
\draw[thick,rounded corners=6pt,gray,thin]
(0cm,.5cm) -- (0cm,4.2cm) -- (5.5cm,4.2cm) -- (5.5cm,0cm) -- (0cm,0cm) -- (0cm,.5cm);
\draw [-arcs'] (0cm,1.9cm) -- (0cm,1.85cm);
\draw [-arcs'] (5.5cm,1.9cm) -- (5.5cm,1.85cm);
\draw [-arcs'] (2.7cm,0cm) -- (2.75cm,0cm);
\draw [-arcs'] (2.7cm,4.2cm) -- (2.75cm,4.2cm);
\end{tikzpicture}}
\put(245,155){\large$B$}
} 
\end{picture}
\end{equation}
}

\noindent According to (\ref{eq:degreecomposition}) we have 
\begin{equation}\label{eq:compAi}
(Id_{X}\times\iota_{D})=(Id_{X}\times\mathcal{A}_{D})\circ(Id_{X}\times\iota)
\end{equation}
and the diagram (\ref{eq:z23}) can be subdivided as in (\ref{eq:z267aart}) below
\\
\begin{equation}\label{eq:z267aart}
\begin{picture}(390,120)(0,80) 

\put(10,173){$X\times S^{k}Y$}
\put(65,175){\vector(1,0){55}}
\put(67,183){$Id_{X}\times\iota$}
\put(130,173){$X\times S^{pk}X$}
\put(190,175){\vector(1,0){50}}
\put(190,183){$Id_{X}\times\mathcal{A}_{D}$}  
\put(250,173){$X\times S^{i}X$}
\put(30,168){\vector(0,-1){55}}
\put(5,137){$\rho_{S^{k}Y}$}
\put(148,170){\vector(0,-1){58}}
\put(152,133){$\mathsf{p}_{S^{pk}X}$}
\put(270,170){\vector(0,-1){60}}
\put(280,133){$\pi_{S^{i}X}$}
\put(19,100){$S^{k}Y$}
\put(50,103){\vector(1,0){70 }}
\put(75,90){$\iota$}
\put(135,100){$ S^{pk}X$}
\put(170,103){\vector(1,0){85}}
\put(205,90){$\mathcal{A}_{D}$}  
\put(260,100){$S^{i}X$}
\put(140,210){$\Delta_{pk}$}
\put(150,180){\begin{tikzpicture}
\pgfsetarrows{inyectivo2->}
\pgfsetlinewidth{.1ex}
\pgfpathmoveto{\pgfpointorigin}
\pgfpathlineto{\pgfpoint{0cm}{-.85cm}}
\pgfusepath{stroke}
\end{tikzpicture}}
\end{picture}
\end{equation}
We will also consider the universal divisor $\Delta_{pk}$ of $S^{pk}X$ and the projections\\
$\pi_{X}:X\times S^{i}X\rightarrow X $,\\
$\mathsf{p}_{X}:X\times S^{pk}X\rightarrow X $ and\\
${\rho_{X}}:X\times S^{k}Y\rightarrow X $.\\
The diagram (\ref{eq:z268art}) involves  the projection $\rho_{S^{k} Y}$ which  is  decomposed as ${\pi}_{S^{k}Y}\circ F $, where $F=f\times Id_{S^{k} Y}$ and ${\pi_{S^{k} Y}}$ is the projection $Y\times S^{k}Y\mapsto S^{k}Y$.
We will often use the line bundle 
\\
${\bar{\xi}}=\rho_{X}^{*}\xi$ on 
$X\times S^{k}Y$.
\\
For example we have 
\begin{equation}\label{eq:idIota110422}
( Id_{X}\times\iota_{D})^{*}(\xi^{-1}(2\Delta_{i}))\cong 
\bar{\xi}
^{-1}(2{\LinL{\Delta}}).
\end{equation}
\begin{equation}\label{eq:z268art}
\begin{picture}(388,200) 

\put(30,173){$X$}
\put(35,168){\vector(0,-1){55}}
\put(27,137){$f$}
\put(29,100){$Y$}
\put(120,173){$X\times S^{k}Y$}
\put(115,175){\vector(-1,0){65}}
\put(80,177){$\rho_{X}$}
\put(138,170){\vector(0,-1){58}}
\put(142,133){$F=f\times Id_{S^{k} Y}$}
\put(123,100){$Y\times S^{k}Y$}
\put(115,103){\vector(-1,0){68}}
\put(80,94){$\pi_{Y}$}
\put(140,100){\vector(0,-1){60}}
\put(130,30){$S^{k}Y$}
\put(145,65){$\pi_{S^{k}Y}$}
\put(230,100){$\rho_{S^{k}Y}$}
\put(150,35){
\begin{tikzpicture}[line width=.5pt]
\draw[thick,rounded corners=6pt,gray,thin]
(0.5cm,5.cm) -- (2.5cm,5.0cm) -- (2.5cm,0cm) -- (0cm,0cm);
\draw [-arcs'] (2.5cm,1.9cm) -- (2.5cm,1.85cm);
\end{tikzpicture}}
\put(95,60){$\Delta_{Y}$}
\put(105,70){\begin{tikzpicture}
\pgfsetarrows{inyectivo3->}
\pgfsetlinewidth{.1ex}
\pgfpathmoveto{\pgfpointorigin}
\pgfpathlineto{\pgfpoint{.85cm}{.85cm}}
\pgfusepath{stroke}
\end{tikzpicture}}
\end{picture}
\end{equation}


\begin{l1}\label{pLchfw+w-}Let $\Delta_{Y}$ be the universal divisor of $S^{k}Y$.Consider the line bundles $\bar{\xi}={\rho_{X}}^{*}\xi$ on $X\times S^{k}Y$ and $\mathcal{L}_{s,j,n}={\pi_{Y}}^{*}\lambda_{s ,j,n}$ on $Y\times S^{k}Y$  . Then we have for $\mathcal{l} \geq 0$

\begin{equation}\label{eq:zz32}
R^{\mathcal{l} }\rho_{S^{k}Y*}({\bar{\xi}}^{s}(-n\LinL{\Delta}))\cong
\bigoplus_{j=0}^{p-1}R^{\mathcal{l} }{\pi}_{S^{k}Y*}(\mathcal{L}_{s ,j,n}(-n\Delta_{Y})).
\end{equation}

\noindent In particular the $\nu^{j}$-eigenbundles of 
$R^{\mathcal{l}}\rho_{S^{k}Y*}({\bar{\xi}}^{s}(-n\LinL{\Delta}))$
and $R^{0}_{F*}\{{\bar{\xi}}^{s}(-n\LinL{\Delta})\}$
are given by 
\begin{equation}\label{eq:zz32b}
R^{\mathcal{l}}\rho_{S^{k}Y*}({\bar{\xi}}^{s}(-n\LinL{\Delta}))(\nu^{j})\cong
R^{\mathcal{l}}{\pi}_{S^{k}Y*}(\mathcal{L}_{s ,j,n}(-n\Delta_{Y}))\mbox{ and}
\end{equation}
\begin{equation}
R^{0}_{F*}\{{\bar{\xi}}^{s}(-n\LinL{\Delta})\}(\nu^{j})\cong \mathcal{L}_{s ,j,n}(-n\Delta_{Y}).
\end{equation}
\end{l1}
\begin{proof}

 Since $\rho_{S^{k} Y}={\pi}_{S^{k}Y}\circ F $, ( $F=f\times Id_{S^{k} Y}$) and  $F$ has finite fibres one can write
 \begin{equation}\label{eq:bschm2n}
\begin{array}{lll}
  R^{\mathcal{l}}\rho_{S^{k}Y*}{\bar{\xi}}^{s}(-n\LinL{\Delta}) &\cong &R^{\mathcal{l}}_{{\pi}_{S^{k} Y}*}(R^{0}_{F*}\{{\bar{\xi}}^{s}(-n\LinL{\Delta})\}).\\
\end{array} 
\end{equation} 
 
\noindent Next we shall write $\LinL{\Delta}$ in terms of the universal divisor
 $\Delta_{pk}$ of $S^{pk}X $ and of the divisor $D$ (see formula (\ref{eq:decundiv}) below). We have 
 (by definition 
of $\LinL{\Delta}$ and equation (\ref{eq:compAi})) 
 \begin{equation}\label{eq:lDDouble}
 \LinL{\Delta}=(Id_{X}\times\iota_{D})^{*}\Delta_{i}=(Id_{X}\times\iota)^{*}(Id_{X}\times\mathcal{A}_{D})^{*}(\Delta_{i})
 \end{equation}
 
\noindent and by the universal property of $\Delta_{i}$ applied to the right-hand side square of the diagram (\ref{eq:z267aart}) one has that $(Id_{X}\times\mathcal{A}_{D})^{*}(\Delta_{i})$
 is the relative divisor of degree $i$ inducing  $\mathcal{A}_{D}$. So we have
 
  \[(Id_{X}\times\mathcal{A}_{D})^{*}(\Delta_{i})= \Delta_{pk}+\mathsf{p}_{X}^{*}D\]
 and using this in (\ref{eq:lDDouble}) we get
 
\begin{equation}\label{eq:decundiv}
\LinL{\Delta} = (\rho_{X})^{*}D+(Id_{X}\times\iota)^{*}\Delta_{pk}. \end{equation}
 Now, using (\ref{eq:decundiv}) in (\ref{eq:bschm2n}) and the fact that

   \[(Id_{X}\times\iota)^{*}\Delta_{pk}=F^{*}\Delta_{Y}\]
  we have
   
   \[
   \begin{array}{lll}
    R^{\mathcal{l} }\rho_{S^{k}Y*}{\bar{\xi}}^{s}(-n\LinL{\Delta})&\cong&  R^{\mathcal{l} }_{{\pi}_{S^{k} Y}*}(R^{0}_{F*}\{{\bar{\xi}}^{s}(-n\rho_{X}^{*}D-nF^{*}\Delta_{Y})\})   \\     
& =& R^{\mathcal{l} }_{{\pi}_{S^{k} Y}*}(R^{0}_{F*}\{{\bar{\xi}}^{s}(-n\rho_{X}^{*}D)\otimes 
\mathcal{O}_{X\times S^{k}Y}(-nF^{*}\Delta_{Y})\}).   \\
\end{array}\]
Applying the projection formula to the direct image $R^{0}_{F*}$ one gets
   \[
   \begin{array}{lll}
R^{\mathcal{l} }\rho_{S^{k}Y*}{\bar{\xi}}^{s}(-n\LinL{\Delta})    &\cong& R^{\mathcal{l} }_{{\pi}_{S^{k} Y}*}(\mathcal{O}_{Y\times S^{k}Y}(-n\Delta_{Y})\otimes R^{0}_{F*}\{ \rho_{X}^{*}(\xi^{s}(-nD)\}).  \end{array}\]
     Now we  use the following base change isomorphism induced from the left-hand side square of (\ref{eq:z268art})      
        \begin{equation}
      R^{0}_{F*}\{ \rho_{X}^{*}(\xi^{s}(-nD))\} \cong  \pi_{Y}^{*}  f_{*}(\xi^{s}(-nD))  
\end{equation}
to get 
 \[\begin{array}{lll}
    R^{\mathcal{l} }\rho_{S^{k}Y*}{\bar{\xi}}^{s}(-n\LinL{\Delta})&\cong & R^{\mathcal{l} }_{{\pi}_{S^{k} Y}*}(\mathcal{O}_{Y\times S^{k}Y}(-n\Delta_{Y})\otimes\pi_{Y}^{*}  f_{*}(\xi^{s}(-nD)) )\\ \end{array} \]

\noindent and using the decomposition into eigenbundles (\ref{eq:z22n})
        we get
          
    \[\begin{array}{lll}
R^{\mathcal{l} }\rho_{S^{k}Y*}{\bar{\xi}}^{s}(-n\LinL{\Delta})
    &\cong&R^{\mathcal{l} }_{{\pi}_{S^{k} Y}*}(\mathcal{O}_{Y\times S^{k}Y}(-n\Delta_{Y})\otimes\pi_{Y}^{*}(\bigoplus^{p-1}_{j=0}\lambda_{s,j,n}))\\
    & =& R^{\mathcal{l} }_{{\pi}_{S^{k} Y}*}(\bigoplus^{p-1}_{j=0}\mathcal{O}_{Y\times S^{k}Y}(-n\Delta_{Y})\otimes 
    \mathscr{L}_{s,j,n})\\
    & =& R^{\mathcal{l} }_{{\pi}_{S^{k} Y}*}(\bigoplus^{p-1}_{j=0}\mathscr{L}_{s,j,n}(-n\Delta_{Y})).
    \end{array} \]
Notice that given $f:X\rightarrow Y$ with $X$ Noetherian, the higher direct images commute with direct sums namely 
$R^{j}f_{*}(\bigoplus_{i} \mathscr{F}_{i})=\bigoplus_{i}R^{j}f_{*}(\mathscr{F}_{i})$.
So we have 
    \begin{equation}
    R^{\mathcal{l} }_{{\pi}_{S^{k} Y}*}(\bigoplus^{p-1}_{j=0}\mathscr{L}_{s,j,n}(-n\Delta_{Y}))=
       \bigoplus^{p-1}_{j=0}R^{\mathcal{l} }_
       {{\pi}_{S^{k} Y}*}(\mathscr{L}_{s,j,n}(- n\Delta_{Y})),
     \end{equation}
from which the Lemma follows.
\end{proof}
\begin{l1}\label{thadd7.4}
Let $Y$ be an irreducible non singular projective curve of genus $g_{Y}$.
Consider the projection $\pi_{S^{k}Y}:Y\times S^{k}Y\rightarrow S^{k}Y$ in (\ref{eq:z268art}) and the universal divisor 
$\Delta_{Y}\subset Y\times S^{k}Y$. For any line bundle 
$M$ on $Y$ and any $m\in \mathds{Z}$ we have the following Chern character
\begin{equation}\label{eq:formula7.4Th}
ch \ \big((\pi_{S^{k}Y})_{!}\pi_{Y}^{*}M(m\Delta_{Y})\big )=((deg \ M+mk+1-g_{Y})-m^{2}\theta)e^{m x}, 
\end{equation}
where $\pi_{Y}^{*}M(m\Delta_{Y})=(\pi_{Y}^{*}M)\otimes 
\mathcal{O}_{Y\times S^{k}Y}(m\Delta_{Y})$.
\end{l1}
\begin{proof} 
This is (7.4)  in \cite{MR1273268} setting $X=Y$, 
$g=g_{Y}$, $i=k$, $\pi=\pi_{S^{k}Y}$ (for Thaddeus 
$\pi:X\times S^{i}X\rightarrow S^{i}X$), 
$k=m\in \mathds{Z}$ and also recall that for Thaddeus 
$\sigma=\theta$ and $\eta=x$. \end{proof}

\begin{proof}[Proof of Theorem \ref{chfw+w-}]
\textbf{a)}  We first notice that 
\begin{equation} 
\begin{array}{lll}\label{eq:bschm}
\iota_{D}^{*}W^{+}_{i}&\cong& R^{1}_{\rho_{S^{k} Y}*}( Id_{X}\times\iota_{D})^{*}(\xi^{-1}(2\Delta_{i})).\\
\end{array}
\end{equation}\\%
By definition (equation (\ref{eq:wimas})) we have 
$\iota_{D}^{*}W^{+}_{i}=\iota_{D}^{*}(R^1_{\pi_{S^{i}X}})_{*}\xi^{-1}(2\Delta_{i})$ then (\ref{eq:bschm}) follows because under our conditions the natural base change morphism 
\begin{equation}\label{eq:bchisomay6}
\iota_{D}^{*}(R^j_{\pi_{S^{i}X}})_{*}\xi^{-1}(2\Delta_{i}) \mapsto R^{j}_{\rho_{S^{k} Y}*}( Id_{X}\times\iota_{D})^{*}(\xi^{-1}(2\Delta_{i})),
\end{equation}
induced from diagram (\ref{eq:z23}), is an isomorphism $\forall \ j\geq 0$. 
\nothing{
Namely, we have that the diagram (\ref{eq:z23}) is cartesian, the sheaf  $\xi^{-1}(2\Delta_{i})$ is flat over $S^{i}X$, the projection $\pi_{S^{i}X}$ is proper and  $S^{i}X$ is locally noetherian.
}
 That is,  from Corollary 2 pg. 50 in \cite{MR0282985}  we see that    the higher direct images $(\!R^j_{\pi_{S^{i}X}}\!)_{*}\xi^{-1}(\!2\Delta_{i}\!)$,  $j\geq 0$, are locally free sheaves and  that for any  $y\in S^{i}X $  the natural maps 
\[\phi^{j}_{y}:(R^j_{\pi_{S^{i}X}})_{*}\xi^{-1}(2\Delta_{i})\otimes k(y)\mapsto H^{j}((X\times S^{i}X)_{y}, \xi^{-1}(2\Delta_{i})_{y})\]
are isomorphisms (in fact, the only ones that are non-zero are  $(\!R^1_{\pi_{S^{i}X}}\!)_{*}\xi^{-1}(\!2\Delta_{i}\!)$ and $\phi^{1}_{y}$).
So the isomorphism in ( \ref{eq:bchisomay6}) follows  as a particular case of Co\-ro\-lla\-ry 6.9.9.2 in \cite{MR163911} or
from Theorem 2.1  in \cite{Conrad-notes}.

\vspace{0.5cm}
\noindent Since (by equation (\ref{eq:idIota110422}))  
$$R^{1}_{\rho_{S^{k} Y}*}( Id_{X}\times\iota_{D})^{*}(\xi^{-1}(2\Delta_{i}))\cong 
R^{1}_{\rho_{S^{k} Y}*}({\bar{\xi}}^{-1}(2\LinL{\Delta})),$$ we have for the $\nu^{j}$-eigenbundles 
\begin{equation}\label{eq:2022-40}
R^{1}_{\rho_{S^{k} Y}*}( Id_{X}\times\iota_{D})^{*}(\xi^{-1}(2\Delta_{i}))(\nu^{j})\cong 
R^{1}_{\rho_{S^{k} Y}*}({\bar{\xi}}^{-1}(2\LinL{\Delta}))(\nu^{j}).
\end{equation}
Now by equation (\ref{eq:bschm}) the left-hand side of 
(\ref{eq:2022-40}) is 
$\iota^{*}_{D}W_{i}^{+}(\nu^{j})$ and the right-hand side is $R^{1}_{{\pi}_{S^{k} Y}*}(\mathcal{L}_{-1,j,-2}(2\Delta_{Y}))$ by 
Lemma \ref{pLchfw+w-} (equation (\ref{eq:zz32b}) taking n=-2, s=-1) . 
 That is
\begin{equation} \iota_{D}^{*}W^{+}_{i}(\nu^{j})\cong R^{1}_{{\pi}_{S^{k} Y}*}(\mathcal{L}_{-1,j,-2}(2\Delta_{Y})),
\end{equation}
 so that using  the Grothendieck-Riemann-Roch Theorem  one can  compute the Chern characters $ch(\iota_{D}^{*}W^{+}(\nu^{j}))$.  In fact,  $ch(\iota_{D}^{*}W^{+}(\nu^{j}))$ can be derived from  Lemma \ref{thadd7.4}, that is, consider 
 $\mathscr{L}_{s,j,n}(2\Delta_{Y})=
\pi^{*}_{Y}\lambda_{s,j,n}\otimes \mathcal{O}
_{Y \times S^{k}Y}(2\Delta_{Y})$ so taking 
$M=\lambda_{s,j,n}$ we have $deg \ M=deg \ \lambda_{s,j,n}$ and in our case $s=-1$ and $n=-2$ so $deg \ M=m'_{j,-2}$ 
($m'_{j,n}$ as defined in the statement of Theorem \ref{chfw+w-}) and so taking $m=2$ Lemma \ref{thadd7.4} tell us that
\[
\begin{split}
ch(R^{0}_{{\pi}_{S^{k} Y}*}(\mathcal{L}_{-1,j,-2}(2\Delta_{Y}))-R^{1}_{{\pi}_{S^{k} Y}*}(\mathcal{L}_{-1,j,-2}(2\Delta_{Y})))=
\ \ \ \ \ \ \ \ \ \  \ \ \ \ \ \ \ \ \ \ 
\\
=((m'_{j,-2}+2k+1-g_{Y})-4\theta)e^{2x},
\end{split}
\]
we will see below that $R^{0}_{{\pi}_{S^{k} Y}*}(\mathcal{O}_{Y\times S^{k}Y}(2\Delta_{Y})\otimes\pi_{Y}^{*}\lambda_{-1,j,-2})=0$
so we have that
\begin{equation*}
ch(R^{1}_{{\pi}_{S^{k} Y}*}(\mathcal{L}_{-1,j,-2}(2\Delta_{Y}))))=-e^{2x}(1+m'_{j,-2}-(-2k+g_{Y}+4\theta)),
\end{equation*}  
that is
 \begin{equation}\label{eq:chimas12ap22}
ch(\iota_{D}^{*}W^{+}_{i}(\nu^{j}))=-e^{2x}(1+m'_{j,-2}-(-2k+g_{Y}+4\theta)).
\end{equation}     
Now $R^{0}_{{\pi}_{S^{k} Y}*}(\mathcal{O}_{Y\times S^{k}Y}(2\Delta_{Y})\otimes\pi_{Y}^{*}\lambda_{-1,j,-2})=0$
because if 
\[
0=R^{\mathcal{0}}_{\rho_{S^{k}Y*}}({\bar{\xi}}^{s}(-n\LinL{\Delta}))
\] 
then $R^{0}_{{\pi}_{S^{k} Y}*}(\mathcal{O}_{Y\times S^{k}Y}(2\Delta_{Y})\otimes\pi_{Y}^{*}\lambda_{-1,j,-2})=0$
because by Lemma \ref{pLchfw+w-} we know that
\[
R^{\mathcal{0} }_{\rho_{S^{k}Y*}}({\bar{\xi}}^{s}(-n\LinL{\Delta}))\cong
\bigoplus_{j=0}^{p-1}R^{\mathcal{0} }_{{\pi}_{S^{k}Y*}}(\mathcal{L}_{s ,j,n}(-n\Delta_{Y})).
\]
Now we have (again by equation (\ref{eq:idIota110422})) that
\[
R^{0}_{\rho_{S^{k} Y}*}( Id_{X}\times\iota_{D})^{*}(\xi^{-1}(2\Delta_{i}))\cong 
R^{0}_{\rho_{S^{k} Y}*}({\bar{\xi}}^{-1}(2\LinL{\Delta})),
\]
and as we mention in the proof of (\ref{eq:bchisomay6}) 
\[
R^{0}_{\rho_{S^{k} Y}*}( Id_{X}\times\iota_{D})^{*}(\xi^{-1}(2\Delta_{i}))=0.
\]
On the other hand, using  (\ref{eq:z21exptheta} ) one has the factorization     \\
\begin{equation}\label{eq:fctCW+vj} 
\frac{e^{\frac{4\theta}{1+2x}}}{(1+2x)^{(1+m'_{j,-2}+2k-g_{Y})}}= (1+2x)^{(r''-g_{Y})}\cdot\prod_{i=1}^{g_{Y}}(1+4\sigma_{i}+2x),
\end{equation}
where $r''=-1-m_{j,-2}-2k+g_{Y}$. So (\ref{eq:fctCW+vj}) can be seen as the Chern class $c(L^{\oplus(r''-g_{Y})}\oplus( L\otimes E))$, where $L$ is a line bundle with Chern class $1+2x$ and $E$ is a rank $g_{Y}$ vector bundle with  $c(E)=e^{4\theta}$. Now we see, by properties of the Chern character and (\ref{eq:chimas12ap22}), that
 \[\begin{array}{lcl}
ch(L^{\oplus(r''-g_{Y})}\oplus( L\otimes E))&=&ch(L^{\oplus(r''-g_{Y})})+ ch( L)ch( E)\\
&=&e^{2x}(r''-g_{Y})+e^{2x}(g_{Y}+4\theta)\\
&=& ch(\iota_{D}^{*}W^{+}_{i}(\nu^{j})).\\
\end{array} 
\]
So we have
 $c(L^{\oplus(r''-g_{Y})}\oplus( L\otimes E))= c(\iota_{D}^{*}W^{+}_{i}(\nu^{j}))$.
 
\vspace{0.5cm}
\noindent
\textbf{proof of b)}: 
\\
Using the definition of $W_{i}^{-}$ (equation (\ref{eq:wimen})), we have
\[
\iota^{*}_{D}W_{i}^{-}=
\iota^{*}_{D}(R^0_{\pi_{S^{i}X}})_{*}\mathscr{E}
\]
where
\[
\mathscr{E}=\{j'_{*}\mathcal{O}_{\Delta_{i}}\}\otimes
\pi_{X}^{*}(\xi)\otimes\mathcal{O}_{X\times S^{i}X}(-\Delta_{i})
\]
is a flat sheaf on $X\times S^{i}X$ because 
$(R^0_{\pi_{S^{i}X}})_{*}\mathscr{E}=W^{-}_{i}$ is a locally free sheaf on $S^{i}X$.
One can argue similarly as in the proof of a) to get the following base change isomorphism induced by the diagram (\ref{eq:z23})
\[
\iota^{*}_{D}(R^0_{\pi_{S^{i}X}})_{*}\mathscr{E}\cong
(R^{0}_{\rho_{S^{k}Y}})_{*}(Id_{X}\times \iota_{D})^{*}\mathscr{E}
\]
that is
\[
\iota_{D}^{*}W^{-}_{i} \cong 
(R^{0}_{\rho_{S^{k}Y}})_{*}(Id_{X}\times \iota_{D})^{*}\mathscr{E}.
\]
One can see that 
\begin{equation}
(Id_{X}\times \iota_{D})^{*}\mathscr{E} \cong \beta'_{*}\{\mathcal{O}_{\bar{\Delta}}\bar{\xi}(-\LinL{\Delta})\}
\end{equation} 
(by using the base change isomorphism $(Id_{X}\times \iota_{D})^{*}j'_{*}\mathcal{O}_{\Delta_{i}}\cong \beta'_{*}\kappa^{*}\mathcal{O}_{\Delta_{i}}$ induced by the square diagram obtained on the top of the diagram (\ref{eq:z23}) by adding the natural projection $\kappa:\LinL{\Delta}\rightarrow \Delta_{i}$). Then

\begin{equation*}
\begin{array}{lll}
\iota_{D}^{*}W^{-}_{i} &\cong & (R^0_{\rho_{S^{k}Y}})_{*}\beta'_{*}\{\mathcal{O}_{\bar{\Delta}}\bar{\xi}(-\LinL{\Delta})\}\\
& &\\
& & \ \ \ \mbox{using the projection formula the last is}\\
& &\\
&\cong&(R^0_{\rho_{S^{k}Y}})_{*}\{\{\beta'_{*}\mathcal{O}_{\bar{\Delta}}\}\otimes\bar{\xi}(-\LinL{\Delta})\}\\
& &\\
& & \ \ \ \mbox{using} \ \ (R^{0}\rho_{S^{k}Y})_{*}=
(R^{0}\pi_{S^{k}Y})_{*}(R^{0}_{F})_{*} \ \ \ \mbox{the last becomes} \\
\end{array}
\end{equation*}
\begin{equation}\label{eq:wimenropif}
\begin{split}
\ \ \ \ \ \ \ \ \ \ \ \cong \ \ &R^{0}_{{\pi}_{S^{k} Y}*}(R^{0}_{F*} \{\{\beta'_{*}\mathcal{O}_{\bar{\Delta}}\}\otimes\bar{\xi}(-\LinL{\Delta})\}).
\end{split}
\end{equation}

\noindent Consider the exact sequence of sheaves on $X\times S^{k}Y$
\begin{equation}\label{eq:z28}
0\rightarrow \bar{\xi}(-2\LinL{\Delta})\rightarrow
\bar{\xi}(-\LinL{\Delta})\rightarrow
\bar{\xi}\otimes\beta'_{*}\mathcal{O}_{\bar{\Delta}}(-\LinL{\Delta})\rightarrow 0.
\end{equation} 
Since $F$ has finite fibres we get an exact sequence of sheaves on 
$Y\times S^{k}Y$
\begin{equation}\label{eq:F*z28}
0\rightarrow R^{0}_{F*} \bar{\xi}(-2\LinL{\Delta})\rightarrow
R^{0}_{F*} \bar{\xi}(-\LinL{\Delta})\rightarrow
R^{0}_{F*} (\bar{\xi}\otimes\beta'_{*}\mathcal{O}_{\bar{\Delta}}(-\LinL{\Delta}))\rightarrow 0,
\end{equation} 
\normalsize
which induces,  for each $j=0,1\cdots p-1$, an exact sequence of $\nu^{j}$-eigen sheaves 
\footnotesize
\begin{equation}\label{eq:F*z28nuj}
0\rightarrow R^{0}_{F*}  \bar{\xi}(-2\LinL{\Delta})(\nu^j)\rightarrow
R^{0}_{F*} \bar{\xi}(-\LinL{\Delta})(\nu^j)\rightarrow
\{R^{0}_{F*} (\bar{\xi}\otimes\beta'_{*}\mathcal{O}_{\bar{\Delta}}(-\LinL{\Delta}))\}(\nu^j)\rightarrow 0.
\end{equation} 
\normalsize
We continue by making the following  observations.\\
1)Let $$\mathscr{F}=R^{0}_{F*} \{\{\beta'_{*}\mathcal{O}_{\bar{\Delta}}\}\otimes\bar{\xi}(-\LinL{\Delta})\}.$$ 
Notice that 
$\mathscr{F}\cong R^{0}_{F*} (\bar{\xi}\otimes\beta'_{*}\mathcal{O}_{\bar{\Delta}}(-\LinL{\Delta}))$ (by the projection formula)
is the non-zero term on the right-hand side of the exact sequence of (\ref{eq:F*z28}). 
\\
2)We have 
$\mathscr{F}=\bigoplus_{j} \mathscr{F}(\nu^{j})$ and then 
$R^{0}_{\pi_{S^{k}Y}*}\mathscr{F}=\bigoplus_{j}R^{0}_{\pi_{S^{k}Y}*}\{\mathscr{F}(\nu^{j})\}$. Since $R^{0}_{\pi_{S^{k}Y}*}\{\mathscr{F}(\nu^{j})\}$ inherits the action of $h$ from 
$\mathscr{F}(\nu^{j})$ we see that 
$$\{R^{0}_{\pi_{S^{k}Y}*}\mathscr{F}\}(\nu^{j})=R^{0}_{\pi_{S^{k}Y}*}\{\mathscr{F}(\nu^{j})\}.$$
3) From (\ref{eq:wimenropif}) $R^{0}_{\pi_{S^{k}Y}*}\{\mathscr{F}\}\cong\iota_{D}^{*}W^{-}_{i}$ so we get that 
\[
\iota_{D}^{*}W^{-}_{i}(\nu^{j})\cong R^{0}_{\pi_{S^{k}Y}*}\{\mathscr{F}(\nu^{j})\}
\] 
for $j=0,\cdots, p-1$. 
\\
4)Also notice that $\mathscr{F}(\nu^{j})\cong\{R^{0}_{F*} 
(\bar{\xi}\otimes\beta'_{*}\mathcal{O}_{\bar{\Delta}}(-\LinL{\Delta}))\}(\nu^{j})$
is the non-zero term on the right-hand side of the exact sequence of (\ref{eq:F*z28nuj}). \\
5) From (\ref{eq:F*z28nuj}) we have the following identity in the Grothendieck ring $K(Y\times S^kY)$ 
$$\mathscr{F}(\nu^{j})=R^{0}_{F*} \bar{\xi}(-\LinL{\Delta})(\nu^j)-R^{0}_{F*}  \bar{\xi}(-2\LinL{\Delta})(\nu^j).$$
So  applying $\pi_{S^{k}Y!}:K(Y\times S^kY)\rightarrow K(S^kY)$ followed by the Chern character $ch$ one has that
\begin{equation}\label{eq:im0dif1}
\begin{split}
ch(\pi_{{S^{k}Y!}}\{\mathscr{F}(\nu^{j})\})=
\hspace{8cm} \ \ \ \ \ 
&\\
ch(\pi_{S^{k}Y!}\{R^{0}_{F*}[\bar{\xi}(-\LinL{\Delta})(\nu^j)]\})
-ch(\pi_{S^{k}Y!}\{R^{0}_{F*} [\bar{\xi}(-2\LinL{\Delta})(\nu^j)]\}).
\hspace{.5cm}
\end{split}
\end{equation}
6) For $s\geq 1$ we have
\begin{equation*}
\begin{split}
0=&(R^s_{{\rho}_{S^{k}Y}})_{*}\beta'_{*}\{\mathcal{O}_{\bar{\Delta}}\bar{\xi}(-\LinL{\Delta})\} \\
&\mbox{similarly to how we arrive to (\ref{eq:wimenropif}) this is}\\
 \cong &R^{s}_{{\pi}_{S^{k} Y}*}(R^{0}_{F*} \{\{\beta'_{*}\mathcal{O}_{\bar{\Delta}}\}\otimes\bar{\xi}(-\LinL{\Delta})\})\\
 =&R^{s}_{\pi_{S^{k}Y}*}\{\mathscr{F}\}\cong R^{s}_{\pi_{S^{k}Y}*}\{\bigoplus_{j}\mathscr{F}(\nu^j)\}\cong\bigoplus_{j} R^{s}_{\pi_{S^{k}Y}*}\{\mathscr{F}(\nu^j)\},\\
\end{split}
\end{equation*}
and consequently $R^{s}_{\pi_{S^{k}Y}*}\{\mathscr{F}(\nu^j)\}=0$.

\vspace{0.5cm}
\noindent With 3) and 6) in (\ref{eq:im0dif1}) above we see that
\begin{equation}\label{eq:chw-nj}
\begin{array}{lll}
\hspace{-.5cm}
\mathftnt{
ch(\iota_{D}^{*}W^{-}_{i}(\nu^j))}=&
\hspace{-.2cm}
\mathftnt{
ch((\pi_{S^{k}Y})_{!}\{R^{0}_{F*} [\bar{\xi}(-\LinL{\Delta})(\nu^j)]\})-
ch((\pi_{S^{k}Y})_{!}\{R^{0}_{F*}[\bar{\xi}(-2\LinL{\Delta})(\nu^j)]\}),
}
\hspace{.5cm} ~
\\
&
\hspace{-.2cm}
\mbox{and using Lemma \ref{pLchfw+w-}  the last is }\\
\hspace{1cm} ~
&
\hspace{-.7cm}
=
\mathftnt{
ch((\pi_{S^{k}Y})_{!}\{\mathcal{L}_{1 ,j,1}(-\Delta_{Y})\})-
ch((\pi_{S^{k}Y})_{!}\{\mathcal{L}_{1 ,j,2}(-2\Delta_{Y})\}).}
\hspace{.7cm} ~
\end{array}
\end{equation}
Applying Lemma \ref{thadd7.4} to each character in (\ref{eq:chw-nj}) with $M = \lambda_{1,j,1}$
for the first one and $M=\lambda_{1,j,2}$ for the second one 
 we have
\begin{equation}\label{eq:z271}
\begin{split}
ch(i^{*}_{D}W^{-}_{i}(\nu^{j}))&=\\
e^{-x}(1+m_{j,1}&-(k+g_{Y}+\theta))-e^{-2x}(1+m_{j,2}-(2k+g_{Y}+4\theta)).
\end{split}
\end{equation}    

\noindent 7) For the calculation of the Chern class we write
\begin{equation}\label{eq:z271-14-04-2022}
 ch(i^{*}_{D}W^{-}_{i}(\nu^{j}))=
e^{-x}(r-\theta))+e^{-2x}(r'+4\theta),\end{equation}
where $r=1+m_{j,1}-k-g_{Y}$ and $r'=-1-m_{j,2}+2k+g_{Y}$.
 Assume $r,r'\geq g_{Y}$. Then 
\[ ch(i^{*}_{D}W^{-}_{i}(\nu^{j}))=
ch((L_{1}\otimes E_{1})  \oplus (L_{2}\otimes E_{2} )),\]
where    $L_{1}$ and $L_{2}$ are line bundles with Chern classes $1-x$ and $1-2x$ respectively and $E_{1}$ and $E_{2}$ are vector bundles  with  Chern characters $r-\theta$ and $r'+4\theta$ respectively.
So,
\[c(i^{*}_{D}W^{-}_{i}(\nu^{j}))=   c(L_{1}\otimes E_{1})  \cdot c(L_{2}\otimes E_{2} ).   \]
From (\ref{eq:z21theta}) one can assume that the non-zero Chern roots of $E_{1}$ and $E_{2}$ are
$-\sigma_1,\dots, -\sigma_{g_{Y}}$ and $4\sigma_1,\dots, 4\sigma_{g_{Y}}$ respectively, so that
{\small \begin{equation}\label{eq:fctcw-}
 c(i^{*}_{D}W^{-}_{i}(\nu^{j}))=((1-x)^{r-g_{Y}}\prod_{i=1}^{g_{Y}}(1-x-\sigma_{i}))\cdot( (1-2x)^{r'-g_{Y}}\prod_{i=1}^{g_{Y}}(1-2x+4\sigma_{i}) )
 \end{equation}
} 
 and  using  (\ref{eq:z21exptheta}) the last is 
 
\begin{equation}\label{eq:clasWmen13ap22}
 =(1-x)^{r}e^{\frac{-\theta}{1-x}}(1-2x)^{r'}e^{\frac{4\theta}{1-2x}}.
\end{equation}
 In case $r,r'\not\geq g_{Y}$ one argues similarly by choosing integers $s,s'$ so that   
 $r+s,r'+s'\geq g_{Y}$, then adding  $e^{-x}s+e^{2x}s'$ to both sides of (\ref{eq:z271-14-04-2022}) one has
\[ 
ch(i^{*}_{D}W^{-}_{i}(\nu^{j}))+e^{-x}(s)+e^{2x}(s')=
e^{-x}(r+s-\theta))+e^{-2x}(r'+s'+4\theta),
\]
which can be seen as
\[
ch(i^{*}_{D}W^{-}_{i}(\nu^{j})\oplus L_{1}^{\oplus s} \oplus L_{2}^{\oplus s'} )=ch((L_{1}\otimes F_{1})\oplus (L_{2}\otimes F_{2})),
\]
where $F_{1}, F_{2}$ have ranks $r+s$ and $r+s'$ respectively.
So
\[
c(i^{*}_{D}W^{-}_{i}(\nu^{j})\oplus L_{1}^{\oplus s} \oplus L_{2}^{\oplus s'} )=c((L_{1}\otimes F_{1})\oplus (L_{2}\otimes F_{2})),
\]
and   by   (\ref{eq:clasWmen13ap22}) we can compute the right-hand side, namely
\[
c(i^{*}_{D}W^{-}_{i}(\nu^{j})\oplus L_{1}^{\oplus s} \oplus L_{2}^{\oplus s'} )=(1-x)^{r+s}e^{\frac{-\theta}{1-x}}(1-2x)^{r'+s'}e^{\frac{4\theta}{1-2x}},
\]
and we also get
\begin{equation}
\begin{split}
c(i^{*}_{D}W^{-}_{i}(\nu^{j}))=
 \ \ \ \ \ \ \ \ \ \ \ \ \  \ \ \ \ \ \ \ \ \ \ \ \ \ \ \ \ \ \ \ \ \ \ \  \ \ \ \ \ \ \ \ \ \   \ \ \ \ \ \ \ \  \ \ \ \ \ \ \ \ &
\\
=(1-x)^{r+s}e^{\frac{-\theta}{1-x}}(1-2x)^{r'+s'}e^{\frac{4\theta}{1-2x}}(c(L_{1}^{\oplus s}))^{-1}(c((L_{2}^{\oplus s'}))^{-1}
\\
=(1-x)^{r+s}e^{\frac{-\theta}{1-x}}(1-2x)^{r'+s'}e^{\frac{4\theta}{1-2x}}(1-x)^{-s}(1-2x)^{-s'}
\\
=(1-x)^{r}e^{\frac{-\theta}{1-x}}
(1-2x)^{r'}e^{\frac{4\theta}{1-2x}}.
\end{split}
\end{equation}   

\end{proof}

 \section{Stable Characteristic Classes}\label{SCC}  
  Theorem \ref{tt3} below is a generalization of Theorem 2.3 in \cite{MR2254543}  where the case $D=0$ is considered. The proof, which we have omitted here,  can be done using similar arguments to those in the proof of Theorem \ref{chfw+w-}  above.

\begin{t1}\label{tt3}
Let $Z_{D}$ be a k-dimensional component of fixed points of $h$ in $S^{i}X$.
Let  $n_{j}$ and $ n_{j}'$ be the degrees of the line bundles  $\lambda_{0,j,-1}$ and 
$\lambda_{0,j,0}$  in formula (\ref{eq:z22n}) respectively. Then 
\begin{equation}\label{eq:z223pa}
\hspace{-.6cm}
a) \ ch(N_{Z_{D}/S^{i}X}(\nu^{j}))=-(1+n'_{j}-g_{Y})+e^{x}(1+n_{j}+k-g_{Y}-\theta),
\end{equation}
\begin{equation}\label{eq:z223pb}
\hspace{-.6cm}
b) \ c(N_{Z_{D}/S^{i}X}(\nu^{j}))=(1+x)^{1+n_{j}+k-g_{Y}}e^{-\frac{\theta}{1+x}},
\end{equation}
\footnotesize
\begin{equation}\label{eq:z84}
c) \ U_{j}(N_{Z_{D}/S^{i}X}(\nu^{j}))=\left(1-\frac{1}{\nu^{j}}\right)^{A}\left(1-\frac{e^{-x}}{\nu^{j}}\right)^{-A}
exp\left(\frac{\theta e^{-x}}{\nu^{j}-e^{-x}}\right)\cdot \left(\frac{1-e^{-x}/\nu^{j}}{1-\nu^{-j}}\right)^{-n_{j}},
\end{equation}
\normalsize
\begin{equation}\label{eq:z85}
d) \ \prod^{p-1}_{j=1}U_{j}(N_{Z_{D}/S^{i}X}(\nu^{j}))=p^{A}m(e^{-x})^{-A}
e^{\theta q(e^{-x})}\prod^{p-1}_{j=1}\left(\frac{1-e^{-x}/\nu^{j}}{1-\nu^{-j}}\right)^{-n_{j}},
\end{equation}
where
 $A=k+1-g_{Y}$, $m(z)=\sum^{p-1}_{j=0}z^{i}$
  and 
   $\ q(z)=\frac{-zm'(z)}{m(z)}$.
\end{t1}

 ~\\
\noindent In particular we have
 \begin{equation}\label{aival} 
{\textnormal{det}(Id-h_{\mid}{N^{\vee}_{Z_{D}/{S^{i}X}}})}={(1-\nu^{p-1})^{a_{1}}\cdot\cdot\cdot (1-\nu)^{a_{p-1}}}
 \end{equation}
where  
\begin{equation*}
a_{j}= rank(N_{Z_{D}/S^{i}X}(\nu^{p-j}))=(n_{p-j}-n'_{p-j}+k) \mbox{~~~ by part a).}
\end{equation*}  

\section{The generalized Chern character for  \texorpdfstring{$B_{i,m,n}$}{lg}}\label{TgCcfBimn}
 
 Let $Z_{D}$ be a component of fixed points of the automorphism $h$. Let $E$ and $F$ be  $h-$linearized vector bundles on  $Z_{D}$.
 The generalized  Chern character of $E$ is given by
 \[ch_{h}(E)=\sum_{j=0}^{p-1}\nu^{j}ch[E(\nu^{j})],\]
 where $ch[E(\nu^{j})]$ is the Chern character of the eigenbundle $E(\nu^{j})$.
 As in the case of the usual Chern character one has that
 \[ch_{h}(E\otimes F)=ch_{h}(E)ch_{h}(F).\]
 So, from equation (\ref{eq:z58}) we have

\begin{equation}\label{eq:chhBimnaTP}
ch_{h}(\iota^{*}_{D}B_{i,m,n})=ch_{h}(\iota^{*}_{D}L^{m}_{i})ch_{h}(\wedge^{i}\iota^{*}_{D}W^{-}_{i})
ch_{h}(\iota^{*}_{D}S^{q_{i}-(i)}U_{i}).
\end{equation}
The factors on the right-hand side of (\ref{eq:chhBimnaTP}) are given in the following result.

\begin{t1}\label{t:chhBimn} Let $Z_{D}$ be a k-dimensional component of fixed points of $h$ in $S^{i}X$.  Let $g_{Y}$ be the genus of the quotient curve $Y$. Let $d_{k}=i-pk$ be the degree of $D$.   Let $\nu^{l}$ and $\nu^{l'}$ be the eigenvalues corresponding to the action of $h$ on the line-bundles $\iota^{*}_{D}L_{i}$ and $\wedge^{i}\iota^{*}_{D}W^{-}_{i}$ respectively.  \\ 
We have the following\\
\begin{equation}\label{eq:z82} \hspace{-1cm}a) \mbox{   } ch_{h}(\iota^{*}_{D}L^{m}_{i})=\nu^{lm} \cdot e^{m(d-2i)x+2mp\theta}.
\end{equation}
\begin{equation}\label{eq:z83} \hspace{-1cm}b)  \mbox{   }  ch_{h}(\wedge^{i}\iota^{*}_{D}W^{-}_{i})
=\nu^{l'}\cdot e^{(d-3i+1-g_{X})x+3p\theta}.
\end{equation}
c)Let $m_{j,n}$ and $m'_{j,n}$  denote respectively the degrees of the bundles $\lambda_{1,j,n}$ and $\lambda_{-1,j,n}$  in formula (\ref{eq:z22n}) .Then 
\begin{equation}\label{eq:z81}
ch_{h}(\iota_{D}^{*}S^{q_{i}-i}U_{i})=
\end{equation}
\footnotesize
\[
=\underset{t^{q_{i}-i}}{coef}\left[exp\left(\frac{-pt^{p}\theta}{e^{px}-t^{p}} \right)\cdot
\frac{(1-t^{p}e^{-px})^{k+g_{Y}-1}}{(1-t^{p}e^{-2px})^{2g_{Y}-2}}
\cdot \prod_{j=0}^{p-1}\left\{\frac{(1-\nu^{j}te^{-2x})^{m_{j,2}+m'_{p-j,-2}}}
{(1-\nu^{j}te^{-x})^{m_{j,1}}} \right\}
\right].
\]
\normalsize

\end{t1}


\begin{proof}
Parts a) and b) follow from  (7.5) in \cite{MR1273268}  and the restriction rules 
$\iota^{*}_{D}\theta=p\theta$  and  $\iota^{*}_{D}x=x$, where  we are using the same notation for the cohomology classes of the symmetric products $S^{i}X$   and $S^{k}Y\cong Z_{D}$  as introduced in Section  $\ref{sec:TchC} $ above 
(recall that for Thaddeus $\sigma=\theta$ and $\eta=x$).\\%
For c),
let $E$ be a rank $r_E$ vector bundle on $Z_{D}$ and let 
\begin{equation}\label{eq:z76a}
P(E,t):=\sum_{l=0}^{\infty}ch[S^{l}(E)]\cdot {t}^{l}.
\end{equation}
One has (see proof of (7.6) in \cite{MR1273268})
\begin{equation}\label{eq:z77e}
P(E,t)=\prod _{\begin{array}{c}
                             \mbox{\scriptsize Chern roots $\alpha$ of $E,$}\\
                          \end{array}} \frac{1}{1-te^{\alpha}}.\\                          
\end{equation}
Let $F$ be an  $h$-linearized vector bundle on $Z_{D}$  and let
\begin{equation}\label{eq:z74}
Q_{h}(F,t)=\sum_{l=0}^{\infty}ch_{h}(S^{l}F)\cdot t^{l}.
\end{equation}
Since $S^{l}F=\bigoplus_{j=0}^{p-1}(S^{l}F)(\nu^{j})$, a Chern  root $\gamma$ of $S^{l}F$ is a Chern root of $(S^{l}F)(\nu^{j})$ for some $j$, say  $\gamma=\sum_{i=1}^{s}\beta_{i}\alpha_{i}$  where $\beta_{i}\geq 0$, $\sum_{i=1}^{s}\beta_{i}=l$ and $\alpha_{i}$ is a Chern root of $F(\nu^{j_{i}})$ for some integer $j_{i}$. Then 
$\nu^{j}e^{\gamma}=m(\nu^{j_{1}}e^{\alpha_{1}},\dots ,\nu^{j_{s}}e^{\alpha_{s}})$ where $m$ is the  degree $l$ monomial $m(x_{1},\dots,x_{s})=\prod_{i=1}^{s}x_{i}^{\beta_{i}}$. 
So one has that
\begin{equation}\label{eq:z77f}
Q_{h}(F,t)=\prod _{\begin{array}{c}
                             \mbox{\scriptsize Chern roots $\alpha$  of $F(\nu^{j}),$}\\
                             \mbox{\scriptsize $j=0,...,p-1$.}
                          \end{array}} \frac{1}{1-\nu^{j}te^{\alpha}},\\                          
\end{equation}
from which one sees that
\begin{equation}\label{eq:z75a}
Q_{h}(F,t)=\prod_{j=0}^{p-1}P(F(\nu^{j}),\nu^{j}t).
\end{equation}
Now we shall assume that $F=\iota_{D}^{*}U_{i}$ and compute $Q_{h}(F,t)$ using (\ref{eq:z75a}).

 \noindent Using (\ref{eq:z55}) and taking Chern class  we have that 

\begin{equation}\label{eq:crUi-17-04-2022}
c(F(\nu^{j}))=c(\iota_{D}^{*}W^{-}_{i}(\nu^{j}))\cdot c((\iota_{D}^{*}W^{+}_{i})^{\vee}(\nu^{j})).
\end{equation}
From equations (\ref{eq:fctcw-}) and (\ref{eq:fctCW+vj})
one has the following factorizations
\begin{equation}\label{eq:crwi-vj-2022}
\begin{split}
c(\iota_{D}^{*}W^{-}_{i}(\nu^{j}))=  \  \ \ \ \ \ \ \ \ \ \ \ \ \ \ \ \ \ \ \ \ \ \ \  \ \ \ \ \ \ \ \ \ \ \ \ \ \ \ \ \ \ \ \ \ \ \ \ \ \ \ \ \ \ \ \ \ \ \ \ \ \ \ \ \\
=(1-x)^{(r-g_{Y})}\cdot (1-2x)^{(r'-g_{Y})}\prod_{i=1}^{g_{Y}}(1-\sigma_{i}-x)\cdot\prod_{i=1}^{g_{Y}}(1+4\sigma_{i}-2x)\\
=\prod _{\begin{array}{c}
                             \mbox{\scriptsize Chern roots $\alpha$ of $\iota_{D}^{*}W^{-}_{i}(\nu^{j})$}\\
                          \end{array}} (1+{\alpha}).\\
\end{split}
\end{equation}
and
\begin{equation}\label{eq:crw+vj-2022}
\begin{split}
 c((\iota_{D}^{*}W^{+}_{i})^{\vee}(\nu^{j}))=
 (1-2x)^{(r''-g_{Y})}\cdot\prod_{i=1}^{g_{Y}}(1-4\sigma_{i}-2x),\\
=\prod _{\begin{array}{c}
                             \mbox{\scriptsize Chern roots $\alpha$ of $(\iota_{D}^{*}W^{+}_{i})^{\vee}(\nu^{j})$}\\
                          \end{array}} (1+{\alpha}).\\
\end{split}
\end{equation}
where 
\begin{equation}\label{eq:vrrprpp2022}
\begin{split}
r & =1+m_{j,1}-k-g_{Y},\ \ \ \ \ \\ r'&=-1-m_{j,2}+2k+g_{Y} \ \ \ \mbox{  and  } \\
r''& =-1-m'_{p-j,-2}-2k+g_{Y}.
\end{split}
\end{equation}\\
So from (\ref{eq:crUi-17-04-2022}), (\ref{eq:crwi-vj-2022}) and (\ref{eq:crw+vj-2022}) we have the Chern roots of $F(\nu^{j})$ and we use them in (\ref{eq:z77e}) to compute
\[
P(F(\nu^{j}),t)=\left(\frac{1}{1-te^{-x}}\right)^{r-g_{Y}}\cdot
\prod_{i=1}^{g_{Y}}\left(\frac{1}{1-te^{-\sigma_{i}-x}}\right)\times
\]
\footnotesize
\[
\times\left(\frac{1}{1-te^{-2x}}\right)^{r'-g_{Y}}\cdot
\prod_{i=1}^{g_{Y}}\left(\frac{1}{1-te^{4\sigma_{i}-2x}}\right)\cdot \left(\frac{1}{1-te^{-2x}}\right)^{r''-g_{Y}}\cdot
\prod_{i=1}^{g_{Y}}\left(\frac{1}{1-te^{-4\sigma_{i}-2x}}\right).
\]
\normalsize
\nothing{
\[
P(U_{i}(\nu^{j}),t)=(\frac{1}{1-te^{-x}})^{1+m_{j,1}-k-g_{Y}-g_{Y}}\cdot
\prod_{i=1}^{g_{Y}}(\frac{1}{1-te^{-\sigma_{i}-x}})\cdot
\]
\[
(\frac{1}{1-te^{-2x}})^{-1-m_{j,2}+2k+g_{Y}-g_{Y}}\cdot
\prod_{i=1}^{g_{Y}}(\frac{1}{1-te^{4\sigma_{i}-2x}})\cdot
\]
\[
(\frac{1}{1-te^{-2x}})^{-1-m'_{p-j,-2}-2k+g_{Y}-g_{Y}}\cdot
\prod_{i=1}^{g_{Y}}(\frac{1}{1-te^{-4\sigma_{i}-2x}}).
\]
}
\noindent Let $h(z):=\frac{1}{1-te^{-z}}$. Expanding the following around $\sigma_{i}=0$ and using $\sigma_{i}^{2}=0$ one has, as the reader may check, that 
\begin{equation}\label{eq:h12022}
h(\sigma_{i}+x)=\frac{1}{1-te^{-\sigma_{i}-x}}=h(x)\left(1+\sigma_{i}\frac{h'(x)}{h(x)}\right),
\end{equation}
\begin{equation}\label{eq:h22022}
h(-4\sigma_{i}+2x)=\frac{1}{1-te^{4\sigma_{i}-2x}}=h(2x)\left(1-4\sigma_{i}\frac{h'(2x)}{h(2x)}\right)
\end{equation}
and 
\begin{equation}\label{eq:h32022}
h(4\sigma_{i}+2x)=\frac{1}{1-te^{-4\sigma_{i}-2x}}=h(2x)\left(1+4\sigma_{i}\frac{h'(2x)}{h(2x)}\right).
\end{equation}
Now we represent  $P(F(\nu^{j}),t)$ as the product $G_{1}G_{2}G_{3}$ of  3 factors defined below which we will modify using (\ref{eq:h12022}), (\ref{eq:h22022}) and (\ref{eq:h32022}):
\[
\hspace{-.8cm}
G_{1}:=
\left(\frac{1}{1-te^{-x}}\right)^{r-g_{Y }}\cdot
\prod_{i=1}^{g_{Y}}\left(\frac{1}{1-te^{-\sigma_{i}-x}}\right)
\]
\[
=\left(\frac{1}{1-te^{-x}}\right)^{r}\cdot
\prod_{i=1}^{g_{Y}}\left(1+\sigma_{i}\frac{h'(x)}{h(x)}\right),
\]
 using (\ref{eq:z21exptheta}) one has that
 $\prod_{i=1}^{g_{Y}}\left(1+\sigma_{i}\frac{h'(x)}{h(x)}\right)=e^{\left(\theta\frac{h'(x)}{h(x)}\right)}$. Also notice that $\frac{h'(x)}{h(x)}=-\frac{te^{-x}}{1-te^{-x}}=-\frac{t}{e^{x}-t}$ so 
\begin{equation}\label{eq:sF12022}
G_{1}=\left(\frac{1}{1-te^{-x}}\right)^{r }\cdot
e^{\left(\theta\frac{h'(x)}{h(x)}\right)}
=(1-te^{-x})^{-r }\cdot
exp\left(\frac{-t\theta}{e^{x}-t}\right).
\end{equation}
In a similar way, we compute
\[G_{2}:=
\left(\frac{1}{1-te^{-2x}}\right)^{r'-g_{Y}}\cdot
\prod_{i=1}^{g_{Y}}\left(\frac{1}{1-te^{4\sigma_{i}-2x}}\right)
\]
\[
=\left(\frac{1}{1-te^{-2x}}\right)^{r'}\cdot
\prod_{i=1}^{g_{Y}}\left(1-4\sigma_{i}\frac{h'(2x)}{h(2x)}\right)
\]

\begin{equation}\label{eq:sF22022}
=\left(\frac{1}{1-te^{-2x}}\right)^{r'}\cdot
e^{\left(-4\theta\frac{h'(2x)}{h(2x)}\right)}
=(1-te^{-2x})^{-r'}\cdot
exp\left(\frac{4t\theta}{e^{2x}-t}\right)
\end{equation}
and
\[G_{3}:=
\left(\frac{1}{1-te^{-2x}}\right)^{r''-g_{Y}}\cdot
\prod_{i=1}^{g_{Y}}\left(\frac{1}{1-te^{-4\sigma_{i}-2x}}\right)
\]
\[
=\left(\frac{1}{1-te^{-2x}}\right)^{r''}\cdot
\prod_{i=1}^{g_{Y}}\left(1+4\sigma_{i}\frac{h'(2x)}{h(2x)}\right)
\]
\begin{equation}\label{eq:sF32022}
=\left(\frac{1}{1-te^{-2x}}\right)^{r''}\cdot
e^{\left(4\theta\frac{h'(2x)}{h(2x)}\right)}
=(1-te^{-2x})^{-r''}\cdot
exp\left(\frac{-4t\theta}{e^{2x}-t}\right).
\end{equation}
Multiplying (\ref{eq:sF12022}), (\ref{eq:sF22022}) and (\ref{eq:sF32022})  we get,  using the values for $r,r'$  and $r''$  given in (\ref{eq:vrrprpp2022}), that
\begin{equation}\label{eq:z76a2}
\begin{split}
P(F(\nu^{j}),t)=\ \ \ \ \ \ \ \ \ \ \ \ \ \ \ \  \ \ \ \ \ \ \ \ \ \ \ \ \ \ \ \ \ \ \ \ \ \ \ \ \ \ \ \ \ \ \ \ \ \ \  \ \ \ \ \ \ \ \ \ \ \ \ \ \ \ \ \ \   &\\
(1-te^{-x})^{-1-m_{j,1}+k+g_{Y}}\cdot exp\left(\frac{-t\theta}{e^{x}-t}\right)
(1-te^{-2x})^{2+m_{j,2}+m'_{p-j,-2}-2g_{Y}}
\end{split}
\end{equation}

\noindent and replacing $t$ by $\nu^{j}t$

\begin{equation}\label{eq:z76}
\begin{split}
P(F(\nu^{j}),\nu^{j}t)=\hspace{1cm} \ \ \ \ \ \ \ \ \ \ \ \ \ \ \ \  \ \ \ \ \ \ \ \ \ \ \ \ \ \ \ \ \ \ \ \ \ \ \ \ \ \ \ \ \ \ \ \ \ \ \  \ \ \ \ \ \ \ \ \ \ \ \ \ \ \ \ \ \   &\\
(1-\nu^{j}te^{-x})^{-1-m_{j,1}+k+g_{Y}}\cdot 
exp\left(\frac{-\nu^{j}t\theta}{e^{x}-\nu^{j}t}\right)(1-\nu^{j}te^{-2x})^{2+m_{j,2}+m'_{p-j,-2}-2g_{Y}}.&\\
\end{split}
\end{equation}
Now we use (\ref{eq:z76}) in (\ref{eq:z75a}) and recalling that $\nu=e^{2i\pi/p}$ one can verify that
\footnotesize
\begin{equation}\label{eq:z80}
Q_{h}( F,t)=exp\left(\frac{-pt^{p}\theta}{e^{px}-t^{p}} \right)\cdot
\frac{(1-t^{p}e^{-px})^{k+g_{Y}-1}}{(1-t^{p}e^{-2px})^{2g_{Y}-2}}
\cdot \prod_{j=0}^{p-1}\left\{\frac{(1-\nu^{j}te^{-2x})^{m_{j,2}+m'_{p-j,-2}}}
{(1-\nu^{j}te^{-x})^{m_{j,1}}} \right\}.
\end{equation}
\normalsize
Therefore
\footnotesize
\[
ch_{h}(\iota_{D}^{*}S^{l}U_{i})=\underset{t^{l}}{coef}(Q_{h}(F,t))
\]
\[
=\underset{t^{l}}{coef}\left[exp\left(\frac{-pt^{p}\theta}{e^{px}-t^{p}} \right)\cdot
\frac{(1-t^{p}e^{-px})^{k+g_{Y}-1}}{(1-t^{p}e^{-2px})^{2g_{Y}-2}}
\cdot \prod_{j=0}^{p-1}\left\{\frac{(1-\nu^{j}te^{-2x})^{m_{j,2}+m'_{p-j,-2}}}
{(1-\nu^{j}te^{-x})^{m_{j,1}}} \right\}
\right].
\]
\normalsize
In particular  if   $l=q_{i}-i$,   ($q_i$ as in  equation (\ref{eq:z58})), one has
\begin{equation}
ch_{h}(\iota_{D}^{*}S^{q_{i}-(i)}U_{i})=\underset{t^{q_{i}-(i)}}{coef}(Q_{h}(F,t))
\end{equation}
\end{proof}
\normalsize
\section{The involution of a hyperelliptic curve}
Putting all data available to us so far in formula (\ref{eq:ec73}) 
 the contribution of a component $Z_{D}$ of fixed points 
in $S^{i}X$ of an automorphism $h$ of order $p$ to the number $N_{i}(h)$ is given by the following 
(see details in Section  \ref{sec:simplification}):
\footnotesize
\begin{multline}\label{eq:z88}
\hspace{-.5cm}C_{i,Z_{D}}(h)=
\frac{p^{A}\nu^{l'+lm}}{(1-\nu^{p-1})^{a_{1}}\cdot\cdot\cdot (1-\nu)^{a_{p-1}}}
\underset{t^{q_{i}-(i)}}{Coef}\mbox{ }\underset{x=0}{Res}\Bigg\{
[m(e^{-x})]^{-A}\prod_{j=1}^{p-1}\left( \frac{1-\frac{e^{-x}}{\nu^{j}}}{1-\nu^{-j}}\right)^{-n_{j}}
\hspace{-.3cm}
\times\\
\hspace{-.5cm}e^{[d(1+m)-i(3+2m)+1-g_{X}]x}\cdot \frac{(1-t^{p}e^{-px})^{k+g_{Y}-1}}{(1-t^{p}e^{-2px})^{2g_{Y}-2}}\cdot
\prod_{j=0}^{p-1}\left\{\frac{(1-\nu^{j}te^{-2x})^{m_{j,2}+m'_{p-j,-2}}}{(1-\nu^{j}te^{-x})^{m_{j,1}}} \right\}
\times\\
\hspace{-.5cm}\left(\frac{x}{1-e^{-x}} \right)^{k-g_{Y}+1}\cdot \frac{\left(1+x\left[q(e^{-x})+p(3+2m)-\frac{pt^{p}}{e^{px}-t^{p}}+\left(\frac{1}{e^{x}-1}-\frac{1}{x} \right)\right] \right)^{g_{Y}}}{x^{k+1}}dx\Bigg\},
\\
\end{multline}
\normalsize
where $A=k+1-g_{Y},$
\begin{equation*}
a_{j}= rank(N_{Z_{D}/S^{i}X}(\nu^{p-j}))=(n_{p-j}-n'_{p-j}+k).
\end{equation*}  
The constants  $l,l',n_{j},n'_{j},m_{j,2},m_{j,1}, m'_{p-j,-2}$  (for their definitions see Theorem \ref{t:chhBimn}, Theorem \ref{tt3} and Theorem \ref{chfw+w-}) appearing in 
 formula  (\ref{eq:z88}) depend on the particular  situation (the curve $X$, the automorphism $h$, the line bundle $\xi$) and in this section we will compute them for the case where $X$ is a hyperelliptic curve of genus $g=g_{X}$, the automorphism $h$ is the hyperelliptic involution and $\xi=K_{X}^{2}$ (see Lemmas \ref{l09-12-2019-22-49}, \ref{l:actionlprime} and \ref{actionLi} below).\\
For the involution of a hyperelliptic curve the contribution  $C_{i,Z_{D}}(h)$ to the Lefschetz number $N_{i}(h)$ does not depend on $D$ but only on the dimension of $Z_{D}$, the dimension of $S^{i}X$ and the genus $g_{X}$ of $X$.
So we write 
\begin{equation}\label{eq:cikgeqCizD}
C_{i,k,g_{X}}(h)=C_{i,Z_{D}}(h)
\end{equation}
 for  $Z_{D}$ a $k-$dimensional component. 
There are $2g_{X}+2$ fixed points of  $h$  in the curve $X$
and 
there are $\binom{2g_{X}+2}{i-2k}$   $k-$dimensional components $Z_{D}$ of fixed points of $h$ in $S^{i}X$ each one corresponding to a divisor $D$ of degree $i-2k$ supported on $i-2k$ distinct fixed  points of $h$. Notice the maximal dimension of a component of fixed points in $S^{i}X$ is $k_{max}=[i/2]$. If we use (\ref{eq:thlsum}) to compute 
the Lefschetz numbers $N_{i}(h)$ then (\ref{eq:3c2p1}) becomes

\begin{equation}\label{eq:tr07-12-2019}Trace( h_{\mid_{ V_{m,n}}} ) =N_{0}(h)+\sum_{i=1}^{w}\sum_{k=0}^{[i/2]}(-1)^{i}\binom{2g_{X}+2}{i-2k} C_{i,k,g_{X}}(h).
\end{equation}

\noindent Let $f_{*}\xi^{s}(-nD)=\bigoplus_{j=0}^1\lambda_{s,j,n}$.

\noindent In Section \ref{sec:hyperellipgticGenus2} we shall use (\ref{eq:tr07-12-2019}) to compute the Verlinde traces of a hyperelliptic curve of genus $g_{X}=2$ and
 to compute the contributions $ C_{i,k,g_{X}}(h)$ it remains to compute the constants $l,l',n_{j},n'_{j},m_{j,2},m_{j,1}, m'_{p-j,-2}$. 
 From their definitions  in Theorem \ref{chfw+w-} and Theorem \ref{tt3} one has:

$m_{j,1}= deg\lambda_{1,j,1}$,

 $m_{j,2}= deg \lambda_{1,j,2}$,

 $m'_{j,n}= deg \lambda_{-1,j,n}$,

$n_{j}=deg\lambda_{0,j,-1}$,

$n_{j}'= deg \lambda_{0,j,0}$.

\noindent In order to compute the degrees of the line bundles $\lambda_{s,j,n}$, consider the virtual representation

$W= H^{0}(X, \xi^{s}(-nD))-H^{1}(X, \xi^{s}(-nD))$,

\noindent then the virtual dimensions of its eigenspaces are given by
\begin{equation}\label{eq:eE2022-1}
dim W(\nu^1) = \frac{1}{2}[L(h^0,\xi^{s}(-nD))-L(h,\xi^{s}(-nD))]
\end{equation}
and
\begin{equation}\label{eq:eE2022-0}
dim W(\nu^0) = \frac{1}{2}[L(h^0,\xi^{s}(-nD))+L(h,\xi^{s}(-nD))].
\end{equation}

\noindent By Riemann-Roch Theorem we have:
\begin{equation}
L(h^0,\xi^{s}(-nD))=sd-n(i-2k)-g_{X}+1.
\end{equation} 

\noindent Next we use the  Atiyah-Bott  formula (\ref{eq:z274ppr}) to compute
$L(h,\xi^{s}(-nD))$
with $\xi=K_{X}^2$ (see (\ref{eq:abap-ppr}) below). Setting $E=\xi^{s}(-nD)$ we have

\begin{equation*}
L(h,\xi^{s}(-nD))=\sum_{l=1}^{2g_{X}+2}\frac{trz \ h|E_{p_{l}}}{det(Id-h|T^{\vee}_{X,p_{l}})}.
\end{equation*}

\noindent Now $h$ acts as multiplication by -1 on  $T^{\vee}_{X,p_{l}}$, so $det(Id-h|T^{\vee}_{X,p_{l}})=2$ .
For  fixed points $p_{l},p_{j}\in X^{h}$ the action of $h$ in the fibre $\mathcal{O}(p_{j})_{p_{l}}$ is multiplication by $(-1)^{\delta_{j,l}}$, where $\delta_{j,l}$ is the  Kronecker delta. To see this, one first notice that the action of $h$ on the fibre $\mathcal{O}_{X}(p_{l})_{p_{l}}=\mathcal{O}_{p_{l}}(p_{l})$ is multiplication by $-1$ because $\mathcal{O}_{p_{l}}(p_{l})$ can be identified with the normal bundle of the point embedding $p_{l}\hookrightarrow X$ which in this case is ${T_{X}}_{,p_{l}}$ since $p_l$ has dimension  zero.
 If $j\not=l$ then when one considers the stalks at $p_{l}$  of the following exact sequence of sheaves induced by the  point embedding $p_{j}\hookrightarrow X$
\[0\rightarrow\mathcal{O}_{X}\rightarrow\mathcal{O}_{X}(p_{j})\rightarrow\mathcal{O}_{p_{j}}(p_{j})\rightarrow 0 \]
one gets an isomorphism on the fibres ${\mathcal{O}_{X}}_{p_{l}}\cong{\mathcal{O}_{X}(p_{j})}_{p_{l}}$ and $h$ acts trivially on the fibres ${\mathcal{O}_{X}}_{p_{l}}$
of the trivial line bundle $\mathcal{O}_{X}$.\\
Now one can compute the action of $h$ on $E_{p_{l}}=(\xi^{s}\otimes\mathcal{O}_{X}(-nD))_{p_{l}}$. First one notice that the action on $(\xi^{s})_{p_{l}}$ is trivial because
we are taking $\xi=K_{X}^2=T_{X}^{\vee}\otimes T_{X}^{\vee}$. As we mentioned at the beginning of this section the divisors $D$ are supported at $\deg D$ distinct points in
 $X^{h}$ then if  $p_{l}$ is a point in the support of $D$ the action of $h$ on 
$\mathcal{O}_{X}(D)_{p_{l}}$ is multiplication by $-1$ and therefore the action on $\mathcal{O}_{X}(-nD)_{p_{l}}$ is multiplication by $(-1)^{n}$. 
 If  $p_{l}$ is not a point in the support of $D$ then one has that the action of $h$  on $\mathcal{O}_{X}(-nD)_{p_{l}}$ is multiplication by $(1)^{n}$. 
 Then when we apply Atiyah-Bott  we get

\begin{equation}\label{eq:abap-ppr}
\begin{array}{lcl}
L(h,\xi^{s}(-nD))&=&(-1)^{n}Deg(D)/2+(1)^n(2g_{X}+2-Deg(D))/2\\
                                        &=&(-1)^{n}(i-2k)/2+  (2g_{X}+2-(i-2k))/2.
                                        \end{array}                                        
\end{equation}

\noindent Next we use that

\[H^{\mathcal{l}}(X, \xi^{s}(-nD))(\nu^j)\cong H^{\mathcal{l}}(Y, \lambda_{s,j,n}),\]

\noindent where  $Y$ is the quotient curve $X/<h> =\mathds{P}^{1}$. We have

\noindent that the Euler characteristics of the eigenbundles $\lambda_{s,j,n}$ are given by

\[\chi(Y, \lambda_{s,j,n})=dim W(\nu^j)=deg\lambda_{s,j,n}+1 \]

\noindent that is,

\begin{equation}
deg\mbox{ } \lambda_{s,j,n}=dim W(\nu^j)-1.
\end{equation}
In particular, we have the following
\begin{l1}\label{l09-12-2019-22-49}
Let $h$ be the involution of a hyperelliptic curve of genus $g_{X}$ and let $\xi= K_{X}^{2}$ then

$m_{1,1}= g_{X}-3$,

$m_{0,1}= 2g_{X}-2+2k-i$,

$m_{1,2}=g_{X}-3-i+2k$,

$m_{0,2}=m_{0,1}$,

$m'_{1,-2}=-2k+i-3g_{X}+1$,

$m'_{2,-2}=m'_{0,-2}=-2k+i-2g_{X}+2$,

$n_{1}=i-2k-g_{X}-1$,

$n'_{1}=-g_{X}-1$.
\end{l1}
\noindent
In the next two lemmas we will use the composition $\iota_{D}$ of equation (\ref{eq:degreecomposition}). 
\begin{l1}\label{l:actionlprime}The action of h on $\wedge^{i}\iota_{D}^{*}W_{i}^{-}$ is multiplication by $(-1)^{i+k}$.
  \end{l1}
\begin{proof} Consider the decomposition into eigenbundles\[\iota_{D}^{*}W_{i}^{-}=  \iota_{D}^{*}W_{i}^{-}(\nu^0)\oplus \iota_{D}^{*}W_{i}^{-}(\nu^1).\]
\noindent
Let  $d_{0}, d_{1}$ be the ranks of the eigenbundles $\iota_{D}^{*}W_{i}^{-}(\nu^0)$ and $\iota_{D}^{*}W_{i}^{-}(\nu^1)$ respectively. Then

\[\wedge^{i}\iota_{D}^{*}W_{i}^{-}=  \wedge^{d_{0}}\iota_{D}^{*}W_{i}^{-}(\nu^0)\otimes \wedge^{d_{1}}\iota_{D}^{*}W_{i}^{-}(\nu^1)\]
\noindent
and the action of the involution $h$ on $\wedge^{i}\iota_{D}^{*}W_{i}^{-}$ is given by

\[\nu^{d_{1}},\]
\noindent
that is,

\[\wedge^{i}\iota^{*}_{D}W_{i}^{-}=  \wedge^{i}\iota^{*}_{D}W_{i}^{-}(\nu^{d_{1}}).\]
\noindent
To compute $d_{1}$ it is enough to compute  degree $0$ part  of $ch(\iota^{*}_{D}W_{i}^{-}(\nu^1))$ (see the expansion of $ch(\mathcal{E})$ in \cite{MR0463157} pg. 432).
So from Theorem \ref{chfw+w-} part b) we have
\begin{equation}
d_{1}=m_{1,1}-m_{1,2}+k,
\end{equation}
then using Lemma \ref{l09-12-2019-22-49} the result follows. 
\end{proof}
\nothing{
the dimension of the fibre of $\iota^{*}W_{i}^{-}(\nu^1)$ at a point  $p\in S^{k}Y$.
As in the proof of Theorem \ref{chfw+w-} b) this can be done by replacing $ch()$ for the dimension of the fibres at $p$ of the respective virtual sheaves in formula (\ref{eq:chw-nj}), that is

\begin{equation}\label{eq:chw-njfiberp}
\begin{array}{lll}
dim(\iota_{D}^{*}W^{-}_{i}(\nu^j))_{p}&=&dim(\rho'_{S^{k}Y}!(\mathcal{L}_{1 ,j,1}(-\Delta_{Y})))_{p}-
dim(\rho'_{S^{k}Y}!(\mathcal{L}_{1 ,j,2}(-2\Delta_{Y})))_{p}\\
&=& \chi (Y,\lambda_{1 ,j,1}(-p))-
\chi (Y,\lambda_{1 ,j,2}(-2p)).\\
\end{array}
\end{equation}

 The dimensions of the eigen spaces $\iota^{*}(W_{i}^{-})_{p}(\nu^j)$, $j=0,1$ are given by

\[dim \mbox{ }\iota^{*}(W_{i}^{-})_{p}(\nu^0) =\frac{1}{2}[Trace({h^{0}\mbox{}}_{\mid}\iota^{*}(W_{i}^{-})_{p}-Trace({h^{1}\mbox{}}_{\mid}\iota^{*}(W_{i}^{-})_{p}],\]

\[dim \mbox{ }\iota^{*}(W_{i}^{-})_{p}(\nu^1) =\frac{1}{2}[{Trace(h^{0}\mbox{}}_{\mid}\iota^{*}(W_{i}^{-})_{p}+Trace({h^{1}\mbox{}}{\mid}\iota^{*}(W_{i}^{-})_{p}]\]

 and the trace of $h^l$ on the fibre  $\iota^{*}(W_{i}^{-})_{p}$ is given by

$ {Trace(h^{l}\mbox{}}_{\mid}\iota^{*}(W_{i}^{-})_{p}= L(h^l,X,K_{X}^2(-p))- L(h^l,X,K_{X}^{2}(-2p))$.

To see that, one can consider the long exact sequence of direct images $R_{\rho_{S^{k}Y}*}^{i}$ associated to the exact sequence  ref{eq:z28}, use the fact that  $\iota^{*}(W_{i}^{-})\cong$ (ref{eq:}) and compute the fibres at $p$.
}   

\begin{l1}\label{actionLi}
The action of $h$ on $\iota_{D}^{*}L_{i}^{m}$ is multiplication by $(-1)^{mi}$.
\end{l1}
\begin{proof}
Let $p\in S^kY$, it will be enough to compute the action on the fibre $ {(L_{i})}_{\iota_{D}(p)}$.  First one notice that
\[
\begin{split}
\mathftnt{
det ((\pi_{S^{i}X})_{!}} 
\mathftnt{\mathcal{O}}_{\mathtiny{X\!\times\!S^{i}X}}
\mathftnt{(\Delta))_{\iota_{D}(p)}=}
\mathftnt{
det(\{{R^{0}_{\pi_{S^{i}X}*}}\mathcal{O}_
{\mathtiny{X\!\times\!S^{i}X}}(\Delta)\}_{\iota_{D}(p)}
\mbox{--}
\{{R^{1}_{\pi_{S^{i}X}*}}\mathcal{O}_
{\mathtiny{X\!\times\!S^{i}X}}(\Delta)\}_{\iota_{D}(p)})}
\\
\hspace{-4.5cm}
\mathftnt{
=det(H^0(X,l(D))-H^1(X,l(D)))}
\hspace{3.7cm}
\\
\hspace{-3.5cm}
\mathftnt{
=det(H^{0}(X,l(D)))\otimes det(H^{1}(X,l(D)))^{-1}},
\hspace{2.5cm}
\end{split}
\]
where  $l=\mathcal{O}_{X}(\iota(p))=\mathcal{O}_{X}(f^{*}p)$, $\iota_{D}$ and $\iota$  as defined in (\ref{eq:degreecomposition}), in particular one has $l(D)=\mathcal{O}_{X}(\iota(p)+D)$.
\\
If we take the virtual representation
\\
$W= H^0(X,l(D))-H^1(X,l(D))$, then
\\
$det\mbox{ } W= det\mbox(W(\nu^{0}))\otimes \det\mbox(W(\nu^{1}))$
and if we consider the virtual dimensions $d_{i}= dim\mbox{ } 
W(\nu^{i})$ then 
$det \mbox (W(\nu^{i}))=(\nu^{i})^{d_{i}}$. So the action 
of $h$ on $det\mbox{ } W$ is given by $(-1)^{d_{1}}$.
These dimensions $d_{i}$ can be computed as explained before $Lemma$ \ref{l09-12-2019-22-49},  that is,  we use (\ref{eq:eE2022-1})  taking $s=1$, $n=-1$, $\xi=l=\mathcal{O}_{X}(f^{*}p)$. One has
\[
L(h^0,l(D))=\deg(l)+\deg(D)-g_{X}+1
\]
\[
{\hspace{2cm}}=2k+(i-2k)-g_{X}+1.
\]
As for $L(h^1,l(D))$ we can use (\ref{eq:abap-ppr})  because the action on the fibres $l_{p_{j}}$ of the fixed points is trivial(the action on the curve $Y$ is trivial and $l=f^{*}\mathcal{O}_{Y}(p)$ is the pull-back of a line bundle on $Y$) exactly as it happens with $\xi=K_{X}^{2}$. Then we obtain
\[
L(h^1,l(D))=(-1)(i-2k)/2+ [2g_{X}+2-(i-2k)]/2.
\]
  Therefore
\[d_{1}=k+(i-2k-g_{x})\equiv i+g_{X}+k {  }\mod 2.\]
\\
Similarly for $\xi = K_{X}^{2}$:

\[det ((\pi_{S^{i}X})_{!}\xi \mathcal{O}_{X\times S^{i}X}(-\Delta))_{\iota_{D}(p)}= det(H^0(X,\xi l^{-1}(-D))-H^1(X,\xi l^{-1}(-D))),\]
\\
and the action of $h$ on $det ((\pi_{S^{i}X})_{!}\xi \mathcal{O}_{X\times S^{i}X}(-\Delta))_{\iota_{D}(p)}$ is given by
\[(-1)^{g_{X}+k}.\]
\\
So the action of $h$ on $\iota^{*}_{D}L_{i}^{m}$ is given by $(-1)^{m(d_{1}+g_{X}+k)}=(-1)^{mi}$.
\end{proof}

\section{A hyperelliptic curve of genus \texorpdfstring{$g_{X}=2$}{lg}}\label{sec:hyperellipgticGenus2}
Let $X$ be a hyperelliptic curve of genus $g_{X}=2$, let\ $h$ be its hyperelliptic involution and 
take $\xi=K^{2}_{X}$. We have  
the embedding $X\overset{\xi K_{X}}{\hookrightarrow} \mathds{P}^{4}$ and we will see that
\begin{equation}
Trace(h_{\mid_{}}{{ V_{l,l(d/2-1)} }})
=dim \ H^{0}(\mathds{P}^{3},\mathcal{O}(l))
\end{equation}
for each integer $l\geq 0$. Notice from equation (\ref{eq:espaciosVerlinde}) that these are the Verlinde 
traces.
Since $\xi$ has degree $d =4g_{X}-4=4$ we have that $w=1$
(recall that $w=[(d-1)/2]$, see Theorem \ref{dimvmn}), then by (\ref{eq:tr07-12-2019}) we have that (taking $m=l$ and $n=l(d/2-1)=l$)
\begin{equation*}
\begin{split}
\mathftnt{
Trace(h_{\mid_{}}{{ V_{l,l(d/2-1)} }})}
=\mathftnt{Trace(h_{\mid_{}}{{ V_{l,l} }})}
=N_{0}(h)+\sum_{i=1}^{1}\sum_{k=0}^{[i/2]}(-1)^{i}\binom{6}{i-2k}C_{i,k,2}(h)\\	
=N_{0}(h) -\binom{6}{1}C_{1,0,2}(h),	\\
	\end{split}
\end{equation*}
one has $m+n= 2l$ so using (\ref{eq:3c2p2i0}) the last becomes
\begin{equation}
\mbox{      }=\underset{t^{2l}}{coef}\left[\frac{1}{\det(I-t\cdot h_{\mid_{}} H^{0}(X,K_{X}\xi))}\right]-\binom{6}{1}C_{1,0,2}(h).
\end{equation}
To compute $\det(I-t\cdot h_{\mid_{}} H^{0}(X,K_{X}\xi))$ we compute the dimensions of the eigenspaces of $h_{\mid_{}} H^{0}(X,K_{X}\xi)$ (this is similar to the calculation of (\ref{eq:eE2022-1}) and (\ref{eq:eE2022-0}) in the previous section) and we have
\[dim H^{0}(X,K_{X}\xi)(\nu^{1})= 3g_{X}-2=4\] 
and
\[dim H^{0}(X,K_{X}\xi)(\nu^{0})= 2g_{X}-3=1.\]

\noindent So we have
\begin{equation}
N_{0}(h)=\underset{t^{2l}}{Coef}
\left( \frac{1}{(1+t)^{4}(1-t)}\right).
\end{equation}
\\
Now, for $C_{1,0,2}(h)$ we have that $i=1$,
$q_{i}-i=n-1=l-1$ and $k=0$. 
\\
By   (\ref{eq:cikgeqCizD}) and (\ref{eq:z88}) 
\begin{equation}
C_{1,0,2}(h)=(-1)^{l+1}\underset{t^{l-1}}{Coef}\cdot\underset{x=0}{Res}\Bigg\{\!\frac{1}{4}\frac{1}{(1-te^{-x})^{2}}
        \left( \frac{1+e^{-x}}{1-e^{-x}}\right)e^{2lx}\frac{(1-te^{-2x})^{2}}{(1+te^{-2x})^{4}}dx\! 
        \Bigg\},
\end{equation}
denoting $e^{-x}$ by $\lambda$ the  residue above becomes
\begin{equation}\label{eq:1}
\underset{\lambda=1}{Res}\Bigg\{ -\frac{1}{4}\frac{1}{(1-t\lambda)^{2}}
        \left(\frac{1+\lambda}{1-\lambda}\right)\frac{1}{\lambda^{2l+1}}\frac{(1-t\lambda^{2})^{2}}
	{(1+t\lambda^{2})^{4}}d\lambda\Bigg\}.
\end{equation}
Notice that if we denote  the function  inside braces by $F(\lambda)$ it has a pole of order $n=1$ on $\lambda=1$, then 
(\ref{eq:1}) is the coefficient of $(\lambda-1)^{n-1}$ in the Taylor expansion of $\left\{(\lambda-1)^{n}F(\lambda)\right\}$ about $\lambda=1$, that is,
\[
        \lim_{\lambda\rightarrow 1}\Bigg[ \frac{1}{(n-1)!}\frac{d^{n-1}}{d\lambda^{n-1}}
        \left\{(\lambda-1)^{n}F(\lambda)\right\}\Bigg].
\]
This limit is equal to $\frac{1}{2(1+t)^{4}}$, in consequence
\[
C_{1,0,2}(h)=(-1)^{l+1}\underset{t^{l-1}}{Coef}\left(\frac{1}{2(1+t)^{4}}\right)
\]
and
\[
N_{1}(h)= \binom{6}{1}C_{1,0,2}(h)=(-1)^{l+1}3\cdot\underset{t^{l-1}}{Coef}\left(\frac{1}{(1+t)^{4}}\right).
\]
We have accomplished 
\begin{equation*}
\begin{split}
Trace(h_{\mid_{}}{{ V_{l,l} }})&=
\underset{t^{2l}}{Coef}\left( \frac{1}{(1+t)^{4}(1-t)}\right)+
(-1)^{l+2}3\cdot\underset{t^{l-1}}{Coef}\left(\frac{1}{(1+t)^{4}}\right)\\
\mbox{\hspace{2cm}}& =\underset{t^{2l}}{Coef}\left( \frac{({1-3t+3t^{2}-t^{3}}}{(1-t^2)^{4}}\right)+
(-1)^{l+2}3\cdot\underset{t^{l-1}}{Coef}\left(\frac{1}{(1+t)^{4}}\right).\\
\end{split}
\end{equation*}
Using the Hilbert series of the ring $K[x_{0},...,x_{n}]$, namely

\begin{equation}\label{eq:21}
	\sum_{d\geq 0}\binom{n+d}{n}t^{d}=\frac{1}{(1-t)^{n+1}},
\end{equation}

\noindent one obtains then
\[Trace(h_{\mid_{}}{{ V_{l,l} }})=\binom{3+l}{3}.
\]
\newpage
\section{Contribution simplification}\label{sec:simplification}
Here we derive formula (\ref{eq:z88}). Using  formula (\ref{eq:ec73})  and the data obtained in Sections \ref{sec:TchC}--\ref{TgCcfBimn}, the contribution of a $k$-dimensional component $Z_{D}$ of ${(S^{i}X)}^{h} $ corresponding to the divisor $D$ is given by 
\[
C_{i,Z_{D}}(h)=
\int_{S^{k}Y} \Bigg\{\nu^{lm} \cdot e^{m(d-2i)x+2mp\theta}
 \times\]
\[
\nu^{l'}\cdot e^{(d-3i+1-g_{X})x+3p\theta}
\times\]
\footnotesize
\[
\underset{t^{q_{i}-i}}{coef}\left[exp\left(\frac{-pt^{p}\theta}{e^{px}-t^{p}} \right)\cdot
\frac{(1-t^{p}e^{-px})^{k+g_{Y}-1}}{(1-t^{p}e^{-2px})^{2g_{Y}-2}}
\cdot \prod_{j=0}^{p-1}\left\{\frac{(1-\nu^{j}te^{-2x})^{m_{j,2}+m'_{p-j,-2}}}
{(1-\nu^{j}te^{-x})^{m_{j,1}}} \right\}
\right]
\times\]
\normalsize
\[
 p^{A}m(e^{-x})^{-A}
e^{\theta q(e^{-x})}\prod^{p-1}_{j=1}\left(\frac{1-e^{-x}/\nu^{j}}{1-\nu^{-j}}\right)^{-n_{j}}
\times\]
\[
\left(\frac{x}{1-e^{-x}}\right)^{k-g_{Y}+1}exp\left(\frac{\theta}{e^{x}-1}-\frac{\theta}{x}\right)
\times
\]
\[
\frac{1}{{(1-\nu^{p-1})^{a_{1}}\cdot\cdot\cdot (1-\nu)^{a_{p-1}}}}
\Bigg\},\]
where the first  3 lines of the integrand correspond to  the generalized Chern character (Theorem \ref{t:chhBimn}), the fourth line is the product of stable characteristic classes of the eigenbundles of the normal bundle $N_{Z_{D}/{S^{i}X}}$ (Theorem \ref{tt3}, d)), the fifth is the Todd class of the component $Z_{D}\cong S^{k}Y$ (formula (\ref{eq:97})) and the sixth line is the inverse of ${\textnormal{det}(Id-h_{\mid}{N^{\vee}_{Z_{D}/{S^{i}X}}})}.$

Now, one can remove $\theta$ from the expression  for $C_{i,Z_{D}}(h)$  by using the following formula (this is  (7.2) from \cite{MR1273268} )
$$
\int_{S^{k}Y}\alpha(x)exp(\beta(x)\theta)= Res_{x=0}\left\{\frac{\alpha(x)(1+x\beta(x))^{g_{Y}}dx}{x^{k+1}}\right\}.
$$
Let $\alpha(x),\beta(x)$ be  as  defined in (\ref{eq:Ax}) and 
(\ref{eq:Bx}) below. 
\footnotesize
\begin{multline}\label{eq:Ax}
\alpha(x)=\Bigg\{
[m(e^{-x})]^{-A}\prod_{j=1}^{p-1}\left( \frac{1-\frac{e^{-x}}{\nu^{j}}}{1-\nu^{-j}}\right)^{-n_{j}}\times
\\
\times
e^{[d(1+m)-i(3+2m)+1-g_{X}]x}\cdot \frac{(1-t^{p}e^{-px})^{k+g_{Y}-1}}{(1-t^{p}e^{-2px})^{2g_{Y}-2}}\times
\\
\times\prod_{j=0}^{p-1}\left\{\frac{(1-\nu^{j}te^{-2x})^{m_{j,2}+m'_{p-j,-2}}}{(1-\nu^{j}te^{-x})^{m_{j,1}}}\right\}
\left(\frac{x}{1-e^{-x}} \right)^{k-g_{Y}+1}\Bigg\}
\end{multline}
\normalsize
and
\footnotesize
\begin{equation}\label{eq:Bx}
\beta(x)=\left[q(e^{-x})+p(3+2m)-\frac{pt^{p}}{e^{px}-t^{p}}+
\left(\frac{1}{e^{x}-1}-\frac{1}{x} \right)\right].
\end{equation}
\normalsize
It is not hard to verify that
$$\hspace{-1cm}C_{i,Z_{D}}(h)=
\frac{p^{A}\nu^{l'+lm}}{(1-\nu^{p-1})^{a_{1}}\cdot\cdot\cdot (1-\nu)^{a_{p-1}}}\int_{S^{i}X} \underset{t^{q_{i}-i}}{coef}\Bigg[ \ \alpha(x)exp(\beta(x)\theta)\ \Bigg].$$
Then
$$\hspace{-1cm}C_{i,Z_{D}}(h)=\frac{p^{A}\nu^{l'+lm}}{(1-\nu^{p-1})^{a_{1}}\cdot\cdot\cdot (1-\nu)^{a_{p-1}}} \times$$
$$\hspace{2cm}\underset{t^{q_{i}-i}}{coef}\Bigg[Res_{x=0}\left\{\frac{\alpha(x)(1+x\beta(x))^{g_{Y}}dx}{x^{k+1}}\right\}\Bigg]
$$
from which we obtained (\ref{eq:z88}).

\vspace{0.5cm}
\noindent{\bf Acknowledgments.} We would like to thank the referee for his comments and suggestions, especially for pointing out several mistakes.
\baselineskip 1mm
\newpage
\addcontentsline{toc}{chapter}{References} 
\bibliographystyle{IEEEtranSA}
\bibliography{references}
\newpage
\appendix

\noindent\textsc{Israel Moreno Mej\'{\i}a.\\
Instituto de Matem\'aticas, Unidad Oaxaca,\\ 
Universidad Nacional Aut\'onoma de M\'exico. \\
Leon 2 altos, Col. Centro, C.P. 68000\\
Oaxaca de Ju\'arez, Oaxaca, M\'exico}\\
\noindent \textit{E-mail address:} 
\textbf{israel@im.unam.mx}\\
\hspace{3cm} \textbf{imm@dunelm.org.uk}\\
\vspace{2cm}    

\noindent\textsc{Dan Silva L\'opez.\\
Instituto de Matem\'aticas, Unidad Oaxaca,\\ 
Universidad Nacional Aut\'onoma de M\'exico. \\
Leon 2 altos, Col. Centro, C.P. 68000\\
Oaxaca de Ju\'arez, Oaxaca, M\'exico}

\vspace{0.3cm}
\noindent Current Address
\\
\vspace{0.1cm}
\noindent\textsc{Instituto Tecnol\'ogico de Huichapan\\
Departamento de Ingenier\'ia en Gesti\'on Empresarial\\
El Saucillo, Huichapan, C.P. 42411\\ 
Hidalgo, M\'exico\\
}
\noindent \textit{E-mail address:} 
\textbf{dsilva@iteshu.edu.mx}\\
\hspace{3cm} \textbf{silva@matem.unam.mx}
\end{document}